\documentclass[12pt]{article}
\usepackage{latexsym,amsmath,amssymb}
\newif\ifdviwin

\usepackage[latin1]{inputenc}
\usepackage[english]{babel}
\usepackage{indentfirst}
\usepackage[mathscr]{eucal}
\usepackage{amssymb,amsmath,amsfonts}
\usepackage{fancybox,fancyhdr}
\usepackage{graphicx}

\newif\ifdviwin

\dviwintrue

\def\cR{\mathcal{R}}

\def\cH{\mathcal{H}}
\def\cT{\mathcal{T}}
\def\cW{\mathcal{W}}

\let\8=\infty \let\0=\emptyset

\let\tilde=\widetilde

\let\landa=\lambda
\let\alfa=\alpha
\let\k=\kappa

\let\parc=\partial

\def\landa{\lambda}

\def\lap{\Delta}
\def\flecha{\rightarrow}
\def\esiz{\langle}
\def\esde{\rangle}

\def\cte.{\mathop{\rm cte.}\nolimits}

\def\cosh{\mathop{\rm cosh }\nolimits}
\def\tanh{\mathop{\rm tanh }\nolimits}

\def\N{\mathbb{N}}
\def\L{\mathbb{L}}

\def\R{\mathbb{R}}

\def\H{\mathbb{H}}
\def\S{\mathbb{S}}
\def\X{\mathfrak{X}}

 \newtheorem{defi}{Definition}
 \newtheorem{teo}[defi]{Theorem}
 
 \newtheorem{cor}[defi]{Corollary}
 \newtheorem{lem}[defi]{Lemma}
 
 \newtheorem{remark}[defi]{Remark}

 \newenvironment{proof}{\rm \trivlist \item[\hskip \labelsep{\it
      Proof}:]}{\par\nopagebreak \hfill $\Box$ \endtrivlist}

\numberwithin{equation}{section}

\begin{document}
\mbox{}\vspace{0.4cm}\mbox{}

\begin{center}
\rule{15.8cm}{1.5pt}\vspace{0.5cm}

\noindent {\large \bf  \sc Hypersurfaces in $\H^{n+1}$ and conformally invariant
equations:\\[0.3cm] the generalized Christoffel and Nirenberg problems}\\ \vspace{0.5cm} {\large José M.
Espinar$\mbox{}^a$, José A. Gálvez$\mbox{}^b$ and Pablo
Mira$\mbox{}^c$}\\ \vspace{0.3cm} \rule{15.8cm}{1.5pt}
\end{center}
  \vspace{1cm}
$\mbox{}^a$, $\mbox{}^b$  Departamento de Geometría y Topología, Universidad de Granada,
E-18071 Granada, Spain. \\ e-mail: jespinar@ugr.es; jagalvez@ugr.es
\vspace{0.2cm}

\noindent $\mbox{}^c$ Departamento de Matemática Aplicada y Estadística,
Universidad Politécnica de Cartagena, E-30203 Cartagena, Murcia, Spain. \\
e-mail: pablo.mira@upct.es \vspace{0.2cm}

\noindent Date: June 11, 2007 \vspace{0.2cm} \\ Keywords: Christoffel problem, Nirenberg problem, Kazdan-Warner conditions, Schouten tensor, hyperbolic Gauss map, Weingarten hypersurfaces.

\vspace{0.3cm}

 \begin{abstract}
Our first objective in this paper is to give a natural formulation of the Christoffel problem for hypersurfaces in $\H^{n+1}$, by means of the hyperbolic Gauss map and the notion of hyperbolic curvature radii for hypersurfaces. Our second objective is to provide an explicit equivalence of this Christoffel problem with the famous problem of prescribing scalar curvature on $\S^n$ for conformal metrics, posed by Nirenberg and Kazdan-Warner. This construction lets us translate into the hyperbolic setting the known results for the scalar curvature problem, and also provides a hypersurface theory interpretation of such an intrinsic problem from conformal geometry. Our third objective is to place the above result into a more general framework.  Specifically, we will show how the problem of prescribing the hyperbolic Gauss map and a given function of the hyperbolic curvature radii in $\H^{n+1}$ is strongly related to some important problems on conformally invariant PDEs in terms of the Schouten tensor. This provides a bridge between the theory of conformal metrics on $\S^n$ and the theory of hypersurfaces with prescribed hyperbolic Gauss map in $\H^{n+1}$. The fourth objective is to use the above correspondence to prove that for a wide family of Weingarten functionals $\mathcal{W} (\k_1,\dots, \k_n)$, the only compact immersed hypersurfaces in $\H^{n+1}$ on which $\mathcal{W}$ is constant are round spheres.
  \end{abstract}

\vspace{0.5cm}

\section{Introduction}

Some of the most interesting problems in the theory of geometric PDEs come from the following classical question: \emph{given a diffeomorphism $G:\S^n\flecha \S^n$ and a smooth function $F:\S^n\flecha \R$, can one find a (necessarily strictly convex) hypersurface $f:\S^n\flecha \R^{n+1}$ with Gauss map $G$ and with $F$ as a prescribed function of its principal curvatures?} Possibly the oldest particular case of this problem is the famous \emph{Christoffel problem} \cite{Chr}, that prescribes $F$ as the mean of the curvature radii of the hypersurface:
 \begin{equation}\label{cris}
F =\frac{1}{n} \sum_{i=1}^n R_i, \hspace{1cm} R_i:= \frac{1}{\k_i},
 \end{equation}
where $\k_1,\dots, \k_n$ are the principal curvatures of the hypersurface. The Christoffel problem has been classically solved after the works \cite{Fir1,Fir2}.

%Our first objective in this paper is to show that the Christoffel problem admits a natural formulation in the context of hypersurfaces in $\H^{n+1}$. The question of giving this satisfactory formulation is unexpectedly subtle, as pointed out by V. Oliker \cite{Oli1}. Indeed, the unit normal of a hypersurface in a space form takes its values in the unit tangent bundle of the space form, and thus cannot be seen as a Gauss map onto $\S^n$. In addition, it is not \emph{a priori} clear if the inverses of principal curvatures will still serve as curvature radii in $\H^{n+1}$. The lack of an immediate formulation for the Christoffel problem in $\H^{n+1}$

It is very natural to ask for the extension of the Christoffel problem to space forms. But surprisingly, even though many interesting contributions on hypersurfaces with prescribed Weingarten curvatures in space forms have been made in the past, a satisfactory development of the Christoffel problem in $\S^{n+1}$ or $\H^{n+1}$ remains unknown. The reason for that seems to be, as Oliker \cite{Oli} points out, that the classical Gauss map is not available on these spaces. Indeed, the unit normal takes its values in the unit tangent bundle of the space form, that is no longer identified with $\S^n$, what makes it unexpectedly subtle even how to formulate the Christoffel problem in $\S^{n+1}$ or $\H^{n+1}$.

%very natural to ask for the extension of the Christoffel problem to space forms. But surprisingly, even though some interesting contributions to that question have been made in the past (...), a totally satisfactory development of the Christoffel problem in space forms remains unknown. For instance, when one tries to formulate the Christoffel problem for hypersurfaces $M^n\subset \H^{n+1}$, the natural candidate to substitute the Gauss map is the \emph{hyperbolic Gauss map} $G:M^n\flecha \S^n$ (see Section 2). However, the hyperbolic Gauss map is not related to the convexity of $M^n$, and thus it is not natural to prescribe a condition like \eqref{cris} that forces the convexity of $M^n$.

Our first goal in this paper is to show that the Christoffel problem can be naturally formulated in the context of hypersurfaces $M^n\subset \H^{n+1}$ in the hyperbolic space. For that we substitute the Euclidean Gauss map by the \emph{hyperbolic Gauss map} $G:M^n\flecha \S^n$, which is widely accepted among specialists in hyperbolic geometry to be the right analogous to the classical Gauss map. In addition, the inverses of the principal curvatures of $M^n\subset \H^{n+1}$ do not serve anymore as curvature radii in this context. We will overcome this difficulty by introducing the \emph{hyperbolic curvature radii} of $M^n\subset \H^{n+1}$, defined as $\mathcal{R}_i:=1/|1-\k_i|$, and that will be shown to play the role in $\H^{n+1}$ of the Euclidean curvature radii from several different perspectives.

The second objective of the paper is to provide a geometric back-and-forth procedure which shows that the Christoffel problem in $\H^{n+1}$ is \emph{essentially} equivalent to a very famous problem from the theory of geometric PDEs. Namely, the Nirenberg problem (or Kazdan-Warner problem) on prescribing scalar curvature in $\S^n$: \emph{given a map $S:\S^n\flecha \R$, does it exist a metric $g=e^{2\rho} g_0$ conformal to the standard metric $g_0$ of $\S^n$, and whose scalar curvature function is given by $S$?} Equivalently, this problem asks for which functions $S$ on $\S^n$ the non-linear elliptic PDE
 \begin{equation}\label{nir}
 \left\{\def\arraystretch{2} \begin{array}{lll}
-\Delta^{g_0} \rho + 1 = \displaystyle \frac{e^{2\rho}}{2} \, S(x) & \text{ if } & n=2, \\
-\Delta^{g_0} u + \displaystyle\frac{n(n-2)}{4} \, u =  \frac{n-2}{4(n-1)} \, S(x) \, u^{\frac{n+2}{n-2}} & \text{ if } & n>2, \hspace{0.3cm} u^{4/(n-2)}:=e^{2\rho},
\end{array} \right.
  \end{equation}
admits a solution globally defined on $\S^n$.

The Nirenberg-Kazdan-Warner problem (\emph{Nirenberg problem} for short from now on) has received an impressive number of contributions over the last 30 years. As a result, researchers on geometric PDEs have clarified to a large extent which smooth functions on $\S^n$ arise as the scalar curvature functions of conformal metrics. We may cite \cite{AmMa,BaCo,BoEz,Cha,ChYa1,ChYa2,ChLi1,ChLi,ChLn,EsSc,KaWa1,KaWa2,Li1,Li2,Mos} as just a few of these works (see the survey \cite{Li5} for more details). However, a complete characterization of scalar curvature functions on $\S^n$ is still unknown.

The back-and-forth construction that we develop here will let us translate all these results on the Nirenberg problem into results for the Christoffel problem in $\H^{n+1}$. And conversely, our construction also provides a hypersurface theory interpretation of an abstract problem of conformal geometry such as the Nirenberg problem. This is not an immediate fact, since the conformal flatness of solutions to the Nirenberg problem is rarely satisfied by the induced metric of a hypersurface in a model space. Besides, a remarkable consequence of \eqref{nir} is that the Christoffel problem in $\H^{n+1}$ cannot be reduced to a linear PDE, in contrast with the classical problem in $\R^{n+1}$.

The third objective of the paper is to analyze the scope of the previous construction. For that, we will consider the \emph{generalized Christoffel problem} in $\H^{n+1}$, in which we prescribe the hyperbolic Gauss map and a given functional of the hyperbolic curvature radii of a compact surface $M^n\subset \H^{n+1}$, not just its mean. Problems of this nature for convex hypersurfaces in $\R^{n+1}$ have been intensely studied, see \cite{GuMa,GLM,GMZ} and references therein. Again, we will show that this question in $\H^{n+1}$ is tightly linked to an important problem from the theory of geometric PDEs that we describe next.

Recall first of all that on a Riemannian manifold $(M^n,g)$, $n>2$, one has the following decomposition: $${\rm Riem} = W_g + {\rm Sch}_g \odot g,$$ where ${\rm Riem}$ is the Riemann curvature tensor, $W_g$ is the Weyl tensor, $\odot$ is the Kulkarni-Nomizu product, and $${\rm Sch}_g :=\frac{1}{n-2}\left({\rm Ric}_g - \frac{S(g)}{2(n-1)} \, g \right)$$ is the \emph{Schouten tensor}. As the Weyl tensor is conformally invariant, the above decomposition reveals that the Schouten tensor encodes all the information on how curvature varies by a conformal change of metric. For this reason the Schouten tensor is the main object of study in conformal geometry. It is also remarkable that $W_g$ vanishes identically in case $(M^n,g)$ is locally conformally flat, which is the situation of the present paper. The \emph{eigenvalues} of ${\rm Sch_g}$ are defined as the eigenvalues of the symmetric endomorphism $g^{-1} {\rm Sch_g}$ obtained by raising an index to ${\rm Sch_g}$.

With this, we will show that the generalized Christoffel problem is equivalent to the problem of prescribing a functional of the eigenvalues of the Schouten tensor for conformal metrics $g= e^{2\rho} g_0$ on $\S^n$, under the regularity condition that $g-2{\rm Sch_g}$ is positive definite. An equivalent version of this regularity condition appeared in \cite{Sch} in connection with the existence of hypersurfaces in $\H^{n+1}$ with a given horospherical metric, but was not linked there with the Schouten tensor of $g$.

A remarkable consequence of this equivalence is that the theory of locally conformally flat Riemannian manifolds can be identified to a large extent with the local theory of hypersurfaces in $\H^{n+1}$ with prescribed regular hyperbolic Gauss map. This fact reports a new way of applying methods from geometric PDEs to investigate hypersurfaces in $\H^{n+1}$, but is also of great interest in the opposite direction: on the one hand, the above equivalence motivates new problems for conformally invariant PDEs that are very interesting from the viewpoint of hypersurfaces in $\H^{n+1}$, although they do not appear so naturally in conformal geometry. And on the other hand, the hypersurface theory interpretation reveals non-trivial \emph{superposition principles} via which one may obtain new solutions to a geometric PDE starting from a previously known one (in the spirit of Backlund transformations, for instance). An example of this use is given in Theorem \ref{inverse}.

At last, the fourth objective of the paper is to give a more definite application of the general correspondence between compact hypersurfaces in $\H^{n+1}$ with regular hyperbolic Gauss map and conformal metrics on $\S^n$. Specifically, we will prove the existence of a wide family of smooth functions $\cW (x_1,\dots, x_n)$ with the following property: \emph{if $M^n\subset \H^{n+1}$ is an immersed compact hypersurface such that $\cW(\k_1,\dots,\k_n)=1$ holds for its principal curvatures, then $M^n$ is a totally umbilical round sphere}. This fact will be implied by a deep theorem in \cite{LiLi1} and the above correspondence.

The above result constitutes a relevant advance in what refers to sphere theorems for Weingarten hypersurfaces, since: (a) we are not assuming \emph{a priori} that the hypersurfaces are embedded, and (b) the family of Weingarten functionals $\cW (x_1,\dots, x_n)$ for which the result holds is extremely large, i.e. it is not just a specific Weingarten relation. We shall also prove a similar classification theorem for horospheres among Weingarten hypersurfaces in $\H^{n+1}$ with one regular end.

The paper is organized as follows. Section 2 will revise the hyperbolic Gauss map for hypersurfaces $M^n\subset \H^{n+1}$ and its relation with tangent horospheres. We will introduce \emph{horospherical ovaloids} as compact hypersurfaces with regular hyperbolic Gauss map and analyze their properties, especially regarding the \emph{horospherical metric} induced on the hypersurface via its associated space of tangent horospheres. Section 3 analyzes the possible formulations of the Christoffel problem in $\H^{n+1}$. We will justify that this problem must be naturally formulated in the class of horospherical ovaloids, and that this leads to the notion of  \emph{hyperbolic curvature radii} $\cR_i :=1/|1-\k_i|$, in terms of which the Christoffel problem in $\H^{n+1}$ is satisfactorily formulated. These two preliminary sections constitute an important part of the paper, because they discuss in great detail why the problem we are considering here seems to be the most natural formulation of the Christoffel problem in $\H^{n+1}$.

In Section 4 we will prove that the Nirenberg problem on $\S^n$ (modulo dilations) is equivalent to the Christoffel problem in $\H^{n+1}$ (modulo parallel translations). This equivalence is made explicit by means of a representation formula for hypersurfaces in terms of the hyperbolic Gauss map and the \emph{horospherical support function}, and shows how to translate into the Christoffel problem the known results for the Nirenberg one.

In Section 5 we generalize the above result, by showing that the generalized Christoffel problem in $\H^{n+1}$ is equivalent to the problem of prescribing a given function of the eigenvalues of the Schouten tensor for conformal metrics on $\S^n$. Some applications of this relationship are explored, in particular regarding problems involving the elementary symmetric functions. In Section 6 we shall prove the above explained characterization of round spheres and horospheres in $\H^{n+1}$ among a very general class of Weingarten hypersurfaces. In the very end of the paper, we will take the inverse approach. and use hypersurface theory in $\H^{n+1}$ to prove: (1) an inversion formula for the eigenvalues of ${\rm Sch}_g$ on an arbitrary locally conformally flat manifold $(M^n,g)$, and (2) a characterization of constant curvature me\-trics on $\S^2$ by the eigenvalues of the $2$-dimensional analogue of the Schouten tensor

From a conceptual viewpoint, this paper investigates hypersurfaces in $\H^{n+1}$ by
analyzing the local variation of their tangent horospheres, thus bifurcating from
the usual perspective in which the variation of tangent hyperplanes is considered.
This approach appears in other works \cite{Eps1,Eps2,Eps3,Sch,FeRo,GMM,GaMi}, but always
from different perspectives. In particular, topics such as the Nirenberg problem or
the Schouten tensor have not been linked to hypersurface theory in $\H^{n+1}$
beforehand.

\section{The hyperbolic Gauss map}

In this preliminary section we study the hypersurfaces in $\H^{n+1}$ with regular hyperbolic Gauss map in terms of their principal curvatures and their tangent horospheres.

\subsection*{Horospheres and the hyperbolic Gauss map}

Let $\H^{n+1}$ denote the $(n+1)$-dimensional hyperbolic space of constant curvature $-1$, and let $\S_{\8}^n=\parc_{\8} \H^{n+1}$ denote its ideal boundary. In what follows, \emph{horospheres} of $\H^{n+1}$ play a central role. These hypersurfaces are easily visualized in the Poincaré ball model $(\mathbb{B}^{n+1},ds^2)$ for $\H^{n+1}$, where here, as usual, $\mathbb{B}^{n+1} =\{x\in \R^{n+1} : ||x||<1\}$. In this model, horospheres correspond to the $n$-spheres that are tangent at one point to the sphere at infinity $\S_{\8}^n$. In this way, two horospheres are always congruent, and they are at a constant distance if their respective points at infinity agree. In addition, given a point $p\in \S_{\8}^n$, the horospheres having $p$ as its point at infinity provide a foliation of $\H^{n+1}$.

From a hypersurface theory viewpoint, horospheres are the flat totally umbilical hypersurfaces in $\H^{n+1}$, and they are complete and embedded.

All of this suggests that horospheres can be naturally regarded in many ways as \emph{hyperplanes} in the hyperbolic space $\H^{n+1}$, even though they are not totally geodesic.

\begin{defi}[\cite{Eps1,Eps2,Bry}]
Let $\phi:M^n\flecha \H^{n+1}$ denote an immersed oriented hypersurface in $\H^{n+1}$ with unit normal $\eta$. The \emph{hyperbolic Gauss map} $$G:M^n\flecha \S_{\8}^n \equiv \S^n$$ of $\phi$ is defined as follows: for every $p\in M^n$, $G(p)\in \S_{\8}^n$ is the point at infinity of the unique horosphere in $\H^{n+1}$ passing through $\phi(p)$ and whose inner unit normal at $p$ agrees with $\eta(p)$.
\end{defi}

Let us point out here that horospheres are globally convex, what allows us to talk about the \emph{inner orientation} of a horosphere, meaning this simply that the unit normal points at the convex side of the horosphere. With respect to this orientation, the second fundamental form of a horosphere is positive definite. Moreover, it turns out that innerly oriented horospheres are the only hypersurfaces in $\H^{n+1}$ with constant hyperbolic Gauss map.

There is an equivalent definition: the \emph{hyperbolic Gauss map} $G:M^n\flecha \S_{\8}^n \equiv \S^n$ of $M^n$ sends each $p\in M^n$ to the point $G(p)$ at the ideal boundary $\S_{\8}^n$ reached by the unique geodesic $\gamma$ of $\H^{n+1}$ that starts at $\phi(p)$ with initial speed $\eta(p)$.

The hyperbolic Gauss map is the analogous concept in the hyperbolic space to the classical Gauss map for hypersurfaces of $\R^{n+1}$, especially if, as we do here, we ask tangent horospheres to play the role of tangent hyperplanes in the Euclidean theory.  It must however be remarked that the \emph{a priori} chosen orientation for the hypersurface matters for the hyperbolic Gauss map. Indeed, if we change the orientation of $M^n$, then $G$ turns into the \emph{negative hyperbolic Gauss map} $G^- :M^n\flecha \S^n$, whose behavior is totally different to that of $G$.

\subsection*{Regularity of the hyperbolic Gauss map}

We shall work in the Minkowski model of $\H^{n+1}$. For that, consider the Minkowski space $\L^{n+2}$ with canonical coordinates $(x_0,\dots, x_{n+1})$ and the Lorentzian metric $$\esiz ,\esde = -dx_0^2 + \sum_{i=1}^{n+1} dx_i^2.$$ The hyperbolic space is then realized in this model as the hyperquadric $$\H^{n+1} = \{x\in \L^{n+2} : \esiz x,x\esde =-1, \ x_0 >0\}.$$ In the same way, the de Sitter $(n+1)$-space and the light cone are given, respectively, by $$ \S_1^{n+1} =\{x\in \L^{n+2} : \esiz x,x\esde =1\}, \hspace{1cm} \N_+^{n+1} = \{x\in \L^{n+2} : \esiz x,x\esde =0, \ x_0 >0\}.$$ Let $\phi:M^n\flecha \H^{n+1}$ be an immersed oriented hypersurface, and let $\eta :M^n\flecha \S_1^{n+1}$ denote its unit normal. Then we can define a normal map associated to $\phi$ taking values in the light cone as
 \begin{equation}\label{light}
 \psi = \phi + \eta : M^n \flecha \N_+^{n+1}.
 \end{equation}
The map $\psi$ is strongly related to the hyperbolic Gauss map $G:M^n\flecha \S_{\8}^n$ of $\phi$. Indeed, the ideal boundary of $\N_+^{n+1}$ coincides with $\S_{\8}^n$, and can be identified with the projective quotient space $\N_+^{n+1} / \R_+$. So, with all of this, we have $G=[\psi] :M^n\flecha \S_{\8}^n\equiv \N_+^{n+1} / \R_+$.

Moreover, if we write $\psi =(\psi_0,\dots, \psi_{n+1})$, then we may interpret $G$ as the map $G:M^n\flecha \S^n$ given by
 \begin{equation}
 G=\frac{1}{\psi_0} (\psi_1,\dots, \psi_{n+1}).
 \end{equation}
In this way, if we label $e^{\rho} := \psi_0$, we get the useful relation
 \begin{equation}\label{trest}
 \psi = e^{\rho} (1,G) :M^n\flecha \N_+^{n+1}.
 \end{equation}
Observe also that, by differentiating \eqref{trest} it follows that
 \begin{equation}\label{estar}
\esiz d\psi, d\psi \esde = e^{2\rho} \esiz dG,dG\esde_{\S^n}.
\end{equation}
We introduce thus the following terminology, in analogy with the Euclidean setting.
\begin{defi}
The smooth function $$e^{\rho}:M^n\flecha \R$$ will be called the \emph{horospherical support function}, or just the \emph{support function}, of the hypersurface $\phi:M^n\flecha \H^{n+1}$.
\end{defi}

%Let $A_{\eta}$ denote the Weingarten operator of the hypersurface $M^n$ in $\H^{n+1}$ with respect to the unit normal $\eta$. Then, as in the Euclidean case, we may define the \emph{principal curvatures} $\{\k_1,\dots, \k_n\}$ and the \emph{principal directions} $\{e_1,\dots, e_n\}$ of the hypersurface at $p\in M^n$ as the eigenvalues and the eigenvectors of $A_{\eta}$ at $p$, respectively. We remark that a horosphere is characterized by the fact that all its principal curvatures are equal to $1$ with respect to its inner orientation.

Besides, if $\{e_1,\dots, e_n\}$ denotes an orthonormal basis of principal directions of $\phi$ at $p$, and if $\k_1, \dots, \k_n$ are the associated principal curvatures, it is immediate that \begin{equation}\label{lemau1}
\esiz d\psi (e_i),d\psi (e_j)\esde = (1- \k_i)^2 \delta_{ij}.
\end{equation}
Thus we have:
\begin{lem}\label{hc}
Let $\phi:M^n\flecha \H^{n+1}$ be an oriented hypersurface. The following conditions are equivalent at $p\in M^n$.
 \begin{enumerate}
 \item[(i)]
The hyperbolic Gauss map $G$ is a local diffeomorphism.
 \item[(ii)]
The associated light cone map $\psi$ in \eqref{light} is regular.
 \item[(iii)]
All principal curvatures of $M^n$ are $\neq 1$.
 \end{enumerate}
\end{lem}
The regularity of the hyperbolic Gauss map gives rise to a notion of convexity specific of the hyperbolic setting, and weaker than the usual geodesic convexity notion:

\begin{defi}[\cite{Sch}]\label{horocon}
Let $M^n\subset \H^{n+1}$ be an immersed oriented hypersurface, and let $\cH_p$ denote the horosphere in $\H^{n+1}$ that is tangent to $M^n$ at $p$, and whose interior unit normal at $p$ agrees with the one of $M^n$. We will say that $M^n$ is \emph{horospherically convex} at $p$ if there exists a neighborhood $V\subset M^n$ of $p$ so that $V\setminus \{p\}$ does not intersect $\cH_p$, and in addition the distance function of the hypersurface to the horosphere does not vanish up to the second order at $p$ in any direction.
\end{defi}
This definition can be immediately characterized as follows.
\begin{cor}\label{schle}
An oriented hypersurface $M^n\subset \H^{n+1}$ is horospherically convex at $p\in M^n$ if and only if all the principal curvatures of $M^n$ at $p$ verify simultaneously $\k_i(p)< 1$ or $\k_i(p)> 1$.
\end{cor}
In particular, if $M^n$ is horospherically convex at $p$ any of the equivalent conditions in Lemma \ref{hc} holds.

\subsection*{Horospherical ovaloids}

\begin{defi}
A compact immersed hypersurface $\phi:M^n\flecha \H^{n+1}$ will be called a \emph{horospherical ovaloid} of $\H^{n+1}$ if it can be oriented so that it is horospherically convex at every point.
\end{defi}
Equivalently, a compact hypersurface is a horospherical ovaloid if and only if it can be oriented so that its hyperbolic Gauss map is a global diffeomorphism. This equivalence follows directly from Lemma \ref{hc} and Corollary \ref{schle} by a simple topological argument, bearing in mind that every compact hypersurface in $\H^{n+1}$ has a point $p$ at which $|\k_i (p)|>1$ for every $i$. In particular, $M^n$ is diffeomorphic to $\S^n$.

It is also immediate from the existence of this point with $|\k_i (p)|>1$ that every horospherical ovaloid has a unique orientation such that $\k_i <1$ everywhere for every $i=1,\dots, n$. We call this orientation the \emph{canonical orientation} of the horospherical ovaloid. It follows that the hyperbolic Gauss map of a canonically oriented horospherical ovaloid is always a global diffeomorphism. This is not necessarily true anymore for the other possible orientation. Let us also point out that if $p$ is a point of a canonically oriented horospherical ovaloid $M^n\subset \H^{n+1}$, then $M^n$ lies around $p$ in the concave part of the unique horosphere that passes through $p$ and whose interior unit normal at $p$ agrees with the unit normal of $M^n$.

Recall that a compact hypersurface $M^n\subset \H^{n+1}$ is a (strictly convex) \emph{ovaloid} if all its principal curvatures are non-zero and of the same sign. Thus, any ovaloid is a horospherical ovaloid, but the converse is not true.

Horospherical ovaloids in $\H^{n+1}$ seem to be an unexplored topic that is of independent interest as a generalization of the usual geodesic ovaloids in $\H^{n+1}$. Nonetheless, let us point out that a horospherical ovaloid is not necessarily embedded. For instance, take a regular curve $\alfa:[0,1]\flecha \H^2$ with geodesic curvature smaller than $1$ at every point, and such that $\alfa(0)=\alfa(1)$ and, moreover, $\alfa'(0)=-\alfa'(1)$. Then by considering $\H^2$ as a totally geodesic surface of $\H^3$ and after rotating $\alfa$ across the geodesic of $\H^2$ that meets $\alfa$ orthogonally at $\alfa(0)$, we get a surface of revolution in $\H^3$ that is a non-embedded horospherical ovaloid.

This lack of embeddedness shows that one cannot talk in general about the outer orientation of a horospherical ovaloid, and justifies the way we introduced the canonical orientation for them.

Another interesting feature of canonically oriented horospherical ovaloids is its good behavior regarding the parallel flow. As usual, the \emph{parallel flow} of an oriented hypersurface $\phi:M^n\flecha \H^{n+1}$ is defined for every $t\in \R$ as $\phi_t :M^n\flecha \H^{n+1}$,
 \begin{equation}\label{lelas}
 \phi_t (p)= {\rm exp}_{\phi(p)} (t \eta_p) : M^n\flecha \H^{n+1},
 \end{equation}
where ${\rm exp}$ denotes the exponential map of $\H^{n+1}$, and $\eta_p$ is the unit normal of $\phi$ at $p$. It is then easy to check that if $\phi$ is a canonically oriented horospherical ovaloid, then the \emph{forward} flow $\{\phi_t\}_t$, $t\geq 0$, is made up by regular canonically oriented horospherical ovaloids. This is no longer true in general for the \emph{backwards} flow (i.e. $t<0$) due to the possible appearance of wave front singularities of the hypersurfaces.

\subsection*{The horospherical metric}

It will be important for our purposes to associate a natural metric to the space of horospheres in $\H^{n+1}$. This construction has appeared in other works previously, but we reproduce it here in order to put special emphasis on some aspects.

Let $\mathcal{M}$ denote the space of horospheres in $\H^{n+1}$. Let us also fix an
arbitrary point $p\in \H^{n+1}$, that we will regard without loss of generality as
the origin in the Poincaré ball model. Then we can view each horosphere
$\mathcal{H}$ as a pair $(x,t)\in \S^n\times \R$, where $x$ is the point at infinity
of $\mathcal{H}$ and $t$ is the (signed) hyperbolic distance of $\mathcal{H}$ to the
point $p$. Here, $t$ is negative if $p$ is contained in the convex domain bounded by
$\mathcal{H}$. Thus we may identify $\mathcal{M}\equiv \S^n\times \R$.

Let us now construct a natural metric on this space of horospheres. Points of the form $(x,0)$ correspond to horospheres passing through the origin in the Poincaré ball model. It is then natural to endow each of these points with the canonical metric $g_0$ of $\S^n$ evaluated at $x$.

But now, the horosphere $(x,t)$ is a parallel hypersurface of $(x,0)$, and the induced metric in $\H^{n+1}$ of this parallel horosphere is a dilation of the one of $\mathcal{H}\equiv (x,0)$, of factor $e^{2t}$. Thus, the natural metric to define at $(x,t)$ is the dilated metric $e^{2t} g_0$ evaluated at $x$. Consequently, we may view the space of horospheres in $\H^{n+1}$ as the product $\S^n\times \R$ endowed with the natural degenerate metric $$\esiz, \esde_{\8} := e^{2t} g_0.$$ Observe that the vertical rulings of $\S^n\times \R$ are null lines with respect to this degenerate metric.

\begin{defi}
Let $\phi:M^n\flecha \H^{n+1}$ denote an oriented hypersurface in $\H^{n+1}$, and let $\mathcal{H}_{\phi} :M^n\flecha \mathcal{M} \equiv \S^n\times \R$ be its tangent horosphere map at every point. We define the \emph{horospherical metric} $g_{\8}$ of $\phi$ as $$ g_{\8} := \mathcal{H}_{\phi} ^* (\esiz,\esde_{\8}),$$ i.e. as the pullback metric via $\mathcal{H}_{\phi}$ of the degenerate metric $\esiz, \esde_{\8}$.
\end{defi}
It turns out that the horospherical metric is everywhere regular if and only if the hyperbolic Gauss map of the hypersurface is a local diffeomorphism. This is a consequence of Lemma \ref{hc} and the following interpretation of the horospherical metric in the Minkowski model of $\H^{n+1}$, i.e. the model in which we will be working.

In the Minkowski model, horospheres of $\H^{n+1}\subset \L^{n+2}$ are the intersections of affine degenerate hyperplanes of $\L^{n+2}$ with $\H^{n+1}$. A simple calculation shows that horospheres are characterized by the fact that its associated light cone map is constant : $\phi + \eta = v\in \N_+^{n+1}$. Moreover, if we write $v= e^{\rho} (1,x)$, we see that $x\in \S^n$ is the point at infinity of the horosphere, and that parallel horospheres correspond to collinear vectors in $\N_+^{n+1}$. This shows that the space of horospheres in $\H^{n+1}$ is naturally identified with the positive null cone $\N_+^{n+1}$. Thus, it is natural to endow this space with the canonical degenerate metric of the light cone, and it is quite obvious from the above construction that this light cone metric coincides with the degenerate metric $\esiz , \esde_{\8}$ defined above.

Consequently, the horospherical metric on a hypersurface in $\H^{n+1}$ is simply the pullback metric of its associated light cone map. Thus, it is regular if and only if the hyperbolic Gauss map is a local diffeomorphism.

All this construction is clearly reminiscent of the usual identification of the
space of oriented vector hyperplanes in $\R^{n+1}$ with the unit sphere $\S^n$. In
this sense, just as the canonical $\S^n$ metric is used in order to measure
geometric quantities associated to the Euclidean Gauss map of a hypersurface in
$\R^{n+1}$, we will use the horospherical metric for measuring geometrical
quantities with respect to the hyperbolic Gauss map. Let us explain in more detail
this consideration, that was first done by Epstein \cite{Eps3}.

First, observe that the ideal boundary $\S_{\8}^n$ of $\H^{n+1}$ does not carry a geometrically useful metric (although it has a natural conformal structure), so we cannot endow $G$ with a pullback metric from the ideal boundary. Nonetheless, let us observe that for defining the hyperbolic Gauss map $G$ we need to know the exact point $p\in \H^{n+1}$ at which we are working (this does not happen for the Euclidean Gauss map). The additional knowledge of this point is then equivalent to the knowledge of the tangent horosphere to the hypersurface at the point. So, it is natural to use the horospherical metric for measuring lenghts associated to the hyperbolic Gauss map. An alternative justification can be found in \cite{Eps3} in connection with the parallel flow of hypersurfaces.

It is interesting to observe that the horospherical metric has played an important role in several different theories. For instance, it is equivalent to the Kulkarni-Pinkall metric \cite{KuPi} (see \cite{Sch}). It also happens that the area of a Bryant surface in $\H^3$ with respect to the horospherical metric is exactly the total curvature of the induced metric of the surface.

\section{The Christoffel problem in $\H^{n+1}$}

The formulation for the Christoffel problem in $\H^{n+1}$ that seems most reasonable at a first sight is: \emph{given a diffeomorphism $G:\S^n\flecha \S^n$ and a function $F:\S^n\flecha \R$, does it exist a hypersurface $\phi:\S^n\flecha \H^{n+1}$ with hyperbolic Gauss map $G$ and such that \eqref{cris} holds for its principal curvatures $\k_1,\dots, \k_n$?} Disappointingly, this is not a natural problem in $\H^{n+1}$, because the two required hypothesis belong to different contexts. Specifically, a compact hypersurface $M^n\subset \H^{n+1}$ whose hyperbolic Gauss map is a global diffeomorphism is not necessarily convex (it is just horospherically convex at every point), and therefore the functional \eqref{cris} may not be defined at some points of $M^n$. On the other hand, the convexity condition that is required on $M^n$ for defining \eqref{cris} is quite meaningless to the hyperbolic Gauss map. These limitations do not appear in the Euclidean setting, where the Gauss map is a global diffeomorphism exactly when the hypersurface is an ovaloid, which is the precise condition needed to define the functional \eqref{cris}.

Other very natural choice \emph{a priori} for acting as curvature radii for the Christoffel problem in $\H^{n+1}$ are the \emph{contact radii}

\begin{equation}\label{curra}
\varrho_i (p):= \coth^{-1} (\k_i (p)).
\end{equation}
These quantities arise from the following fact: if $\alfa$ is a curve in $\H^2$ with geodesic curvature $k_g$ at $p$, then the inverse of the curvature of the unique circle in $\H^2$ having an order two contact with $\alfa$ at $p$ is given by $\coth^{-1} (k_g)$. Nonetheless, we again see that the quantities \eqref{curra} are not well defined if $\k_i \in [-1,1]$, so we also have to discard them.

The above discussion concludes that we must seek an alternative formulation of the Christoffel problem in $\H^{n+1}$. More specifically, we need to find a more adequate notion of curvature radii for hypersurfaces in $\H^{n+1}$ that makes sense exactly when the hypersurface is a horospherical ovaloid.

\subsection*{The hyperbolic curvature radii}

The key observation at this point is that the Euclidean Gauss map $N:M^n\flecha \S^n$ of a strictly convex hypersurface $M^n\subset \R^{n+1}$ is related in a very simple way to the Euclidean curvature radii $R_i$, as follows. Let $$\alfa_i (t): (-\epsilon,\epsilon)\flecha M^n$$ denote a curve in $M^n\subset \R^{n+1}$ with $\alfa_i (0)=p$ and such that $\alfa_i' (0)= e_i$, where $\{e_1,\dots, e_n\}$ is an orthonormal basis of principal directions in $T_p M^n$. Let $L_0^t (\alfa_i)$ (resp. $L_0^t (N\circ \alfa_i)$) denote the length of $\alfa_i ([0,t])$  (resp. of $N\circ \alfa_i([0,t])$), where $N$ is the unit normal of $M^n$, that is assumed to be a local diffeomorphism at $p$. Then
 \begin{equation}\label{ahiva}
\lim_{t\to 0} \frac{L_0^t (\alfa_i)}{L_0^t (N\circ \alfa_i)} = \lim_{t\to 0} \frac{\int_0^t |\alfa_i ' (u)| du}{\int_0^t |(N\circ \alfa_i) ' (u)| du} = \frac{|\alfa_i'(0)|}{|(N\circ \alfa_i)' (0)|}= \frac{1}{|\k_i(p)|} = R_i (p).
 \end{equation}
This relation is relevant to the Christoffel problem, since it indicates that the curvature radii admit an interpretation in terms of the Gauss map.

The above construction can also be carried out in $\H^{n+1}$, with the Euclidean Gauss map substituted by the hyperbolic Gauss map. Specifically, let us consider an oriented hypersurface $M^n\subset \H^{n+1}$ that is horospherically convex at $p\in M^n$, and let $\{e_1,\dots, e_n\}$ be an orthonormal basis of principal directions in $T_pM^n$. If we now take $\alfa_i :(-\epsilon, \epsilon)\flecha M^n$ a regular curve with $\alfa_i (0)=p$ and $\alfa_i' (0)= e_i$, we may define in analogy with the Euclidean situation the \emph{hyperbolic curvature radii} of $M^n$ at $p$ as $$\cR_i (p):= \lim_{t\to 0} \frac{L_0^t (\alfa_i)}{L_0^t (G\circ \alfa_i)}.$$ Here the length of the hyperbolic Gauss map along $\alfa_i$ is obviously taken with respect to the horospherical metric $g_{\8}$. In addition, the quotients make sense since $G$ is a local diffeomorphism at $p$. At last, by an argument analogous to \eqref{ahiva} we get $\cR_i(p) = 1/|1-\k_i (p)|.$ So, we propose the following definition as the natural analogue in $\H^{n+1}$ of the Euclidean curvature radii for geometrical problems involving the hyperbolic Gauss map.

\begin{defi}
Let $M^n\subset \H^{n+1}$ be a hypersurface that is horospherically convex at $p\in M^n$. We define the \emph{hyperbolic curvature radii} $\{\cR_1,\dots, \cR_n\}$ of $M^n$ at $p$ as $$\cR_i(p) = \frac{1}{|1-\k_i (p)|},$$ where here $\{\k_1,\dots, \k_n\}$ are the principal curvatures of $M^n$ at $p$.
\end{defi}

\begin{remark}
We may observe that $1/\cR_i$ is simply the length of $dG_p (e_i)$ with respect to the horospherical metric, where here $e_i$ is a principal unit vector at $p$. As the same property is true in the Euclidean setting, this indicates again that $\cR_i$ is the proper extension to $\H^{n+1}$ of the Euclidean curvature radii.
\end{remark}

The above arguments suggest the following formulation of the Christoffel problem in $\H^{n+1}$ as the most natural one: \emph{given a diffeomorphism $G:\S^n\flecha \S^n$ and a smooth function $C:\S^n\flecha \R$, does it exist an oriented hypersurface $\phi:\S^n\flecha \H^{n+1}$ (necessarily a horospherical ovaloid) with hyperbolic Gauss map $G$ and with $C$ as the mean of the hyperbolic curvature radii?} It must be stressed here that the hyperbolic curvature radii $\cR_i$ make sense in the compact case exactly for the class of horospherical ovaloids, which is the condition we were looking for.

It is convenient to work with a simplified (but equivalent) version of the above Christoffel problem. In order to expose this simplification (see next page), we introduce the next remark, as well as the following subsection.

\begin{remark}
As $G$ is a global diffeomorphism in the Christoffel problem, it can be used as a global parametrization of the horospherical ovaloid. In other words, we may assume that $G(x)= x$ on $\S^n$ without losing generality.
\end{remark}

\subsection*{Orientation and the parallel flow}

It is interesting to observe the behaviour of the Christoffel problem under the parallel flow. Let $\{\phi_t\}_{t\in \R}$ denote the parallel flow of a solution $\phi:\S^n\flecha \H^{n+1}$ to the Christoffel problem for the function $C(x)$. Then the hyperbolic Gauss map $G(x)=x$ remains invariant under this flow, and the horospherical metric of $\phi_t$ is $g_{\8,t} = e^{2t }g_{\8}$. Moreover, the principal curvatures $\k_i^t$ of $\phi_t$ at regular points are given by
 \begin{equation}\label{curprin}
 \k_i^t (p)=  \frac{\k_i (p) -\tanh (t)}{1-\k_i (p) \tanh (t)},
 \end{equation}
and so the mean $C_t(x)$ of the hyperbolic curvature radii of $\phi_t$ is, at its regular points,
 \begin{equation}\label{cesubte}
 C_t (x)= \frac{1}{2} - \frac{e^{-2t}}{2} \left(1- 2 C(x)\right).
 \end{equation}
Unfortunately, $\phi_t$ is not always regular. Indeed, the first fundamental form of $\phi_t$ is given by
 \begin{equation}\label{1ffpar}
I_t (e_i,e_j) = (\cosh t -\k_i \sinh t)^2 \delta_{ij},
\end{equation}
where $\{e_1,\dots, e_n\}$ is an orthonormal basis of principal directions of $M^n$, and it can be singular. Nonetheless, it is immediate from \eqref{1ffpar} the existence of some $t_0\in \R$ such that $\phi_t$ is regular (and hence solves the Christoffel problem for the function $C_t (x)$ in \eqref{cesubte}) for $t\geq t_0$. In this way, the solutions to the Christoffel problem come in $1$-parameter families determined by the parallel flow in the above way. Moreover, if $\k_i <1$, i.e. the solution to the Christoffel problem is canonically oriented, $\phi_t$ is regular for every $t\geq 0$.

\begin{remark}
In the Euclidean Christoffel problem it is usually assumed that the ovaloid $M^n\subset \R^{n+1}$ is canonically oriented, so that $R_i = 1/\k_i >0$. In contrast, in the hyperbolic case it is at a first sight restrictive to deal only with canonically oriented horospherical ovaloids $M^n$ in $\H^{n+1}$, as a change of orientation on $M^n$ transforms the hyperbolic Gauss map $G$ into the \emph{negative} hyperbolic Gauss map $G^-$, which is totally different from $G$.

Nevertheless, if $\phi:\S^n\flecha \H^{n+1}$ is a negatively oriented horospherical ovaloid, there exists $t_0>0$ such that if $t\geq t_0$ the parallel hypersurface $\phi_t :\S^n\flecha \H^{n+1}$ is a canonically oriented horospherical ovaloid.

This property lets us work without loss of generality with the canonically oriented situation in the Christoffel problem in $\H^{n+1}$, as this problem is invariant under the parallel flow.
\end{remark}

\subsection*{Formulation of the Christoffel problem in $\H^{n+1}$}

Taking all of this into account, we formulate the Christoffel problem as follows:

\vspace{0.3cm}

\noindent{\bf The Christoffel problem in $\H^{n+1}$:}

\begin{quote}
Let $C:\S^n\flecha \R_+$ denote a positive smooth function. Find out if there exists a canonically oriented horospherical ovaloid $\phi:\S^n\flecha \H^{n+1}$ such that its hyperbolic Gauss map and its mean of the hyperbolic curvature radii are, respectively,
 \begin{equation}\label{crisdata}
G(x)=x  \hspace{0.6cm} \text{ and } \hspace{0.6cm} C(x)= \frac{1}{n} \sum_{i=1}^n \frac{1}{1-\k_i}
 \end{equation}
for every $x\in \S^n$. Here $\k_1,\dots, \k_n$ are the principal curvatures of $\phi$.
\end{quote}

\begin{defi}
We will say that a smooth function $C:\S^n\flecha \R_+$ is a \emph{Christoffel function} if it arises as the mean of the hyperbolic curvature radii of some canonically oriented horospherical ovaloid $\phi$ in $\H^{n+1}$ with hyperbolic Gauss map $G(x)=x$.
\end{defi}

\noindent {\bf Convention:} From now on, and unless otherwise stated, by a horospherical ovaloid we will always mean a canonically oriented one. So, we will have
 \begin{equation}\label{crisfun}
 \cR_i =\frac{1}{1-\k_i}.
 \end{equation}
 
In \cite{Oli2} one can find an alternative formulation of a Christoffel problem in $\H^{n+1}$, which is very different from the one here. However, the Christoffel-type problem proposed in that work is not invariant by isometries of $\H^{n+1}$ (since it implicitly uses an \emph{Euclidean Gauss map}), while the one formulated here does not have this disadvantage.

\section{Solution of the Christoffel problem}

This section is devoted to show the equivalence of the Christoffel problem in $\H^{n+1}$ and the Nirenberg problem (or Kazdan-Warner problem) on prescribing scalar curvature on $\S^n$, by means of an explicit back-and-forth procedure.

\subsection*{A representation formula}

Firstly, we shall deduce a formula that represents locally a hypersurface in $\H^{n+1}$ in terms of its hyperbolic Gauss map $G$ and its support function $e^{\rho}$. As we shall work at a point around which $G$ is a local diffeomorphism, we may use the hyperbolic Gauss map to parametrize the hypersurface, i.e. we may assume that the hypersurface is given by $\phi:U\subset \S^n\flecha \H^{n+1}$ with $G(x)=x$, where $U$ is an open set of $\S^n$.

\begin{teo}\label{represe}
Let $\phi:U\subset \S^n\flecha \H^{n+1}$ denote a locally horospherically convex hypersurface with hyperbolic Gauss map $G(x)=x$, and support function $e^{\rho}:U\flecha (0,+\8)$. Then it holds
 \begin{equation}\label{repfor}
 \phi = \frac{e^{\rho}}{2}\left( 1+ e^{-2\rho} \left( 1+ ||\nabla^{g_0} \rho ||_{g_0}^2 \right)\right) (1,x) + e^{-\rho} (0, -x +\nabla^{g_0} \rho).
 \end{equation}
\end{teo}
\begin{proof}
We will give a constructive proof, in order to explain the origin of the expression \eqref{repfor}. Let $g:=e^{2\rho} g_0$, where $g_0$ is the canonical metric in $\S^n$. It follows from \eqref{estar} that $\esiz d\psi,d\psi\esde = e^{2\rho} \esiz dG,dG\esde_{\S^n} = e^{2\rho} g_0=g$. Our first objective is to prove the formula
 \begin{equation}\label{lapla}
 \phi = \frac{1}{n} \lap^g \psi + \left(\frac{S(g) + n(n-1)}{2n (n-1)}\right) \psi,
 \end{equation}
where $S(g)$ is the scalar curvature of $g$, and $\lap^g \psi$ stands for the Laplacian of $\psi$ with respect to $g$.

Let $\{e_1,\dots, e_n\}$ be an orthonormal basis of principal directions at $p$, with principal curvatures $\k_1,\dots, \k_n$. If we label $v_i :=1/(1-\k_i) \, e_i$, it follows that $\esiz d\psi (v_i),d\psi (v_j)\esde = \delta_{ij}$. Let us now view $\psi :M^n\flecha \N_+^{n+1}\subset \L^{n+2}$ as a spacelike codimension-$2$ submanifold of $\L^{n+2}$. Its normal space at every point is spanned by $\{\phi,\eta\}$, and its second fundamental form $\alfa :\X (M^n)\times \X(M^n)\flecha \L^{n+2}$ is given by
 \begin{equation}\label{2ff}
 \alfa (v_i,v_j) = \left( \frac{1}{1-\k_i} \, \phi + \frac{\k_i}{1-\k_i} \, \eta \right) \delta_{ij}.
 \end{equation}
If $K(x,y)$ denotes the sectional curvature of $\psi$, the Gauss equation in $\L^{n+2}$ lets us infer from \eqref{2ff} that
 \begin{equation*}
K(v_i,v_j) =\esiz \alfa (v_i,v_i),\alfa (v_j,v_j)\esde - ||\alfa (v_i,v_j)||^2 = 1- \frac{1}{1-\k_i} - \frac{1}{1-\k_j}.
 \end{equation*}
Hence  \begin{equation}\label{kei}
 S(g) = n(n-1) - 2(n-1) \sum_{i=1}^n \frac{1}{1-\k_i}.
 \end{equation}
Also by \eqref{2ff} we get that the mean curvature vector of $\psi$ in $\L^{n+2}$ is $$\def\arraystretch{1.5}\begin{array}{lll}{\bf H} &=& \displaystyle\frac{1}{n} \sum_{i=1}^n \alfa (v_i,v_i) = \displaystyle \frac{1}{n} \sum_{i=1}^n \left(\displaystyle \frac{1}{1-\k_i} \, \phi + \displaystyle\frac{\k_i}{1-\k_i} \, \eta \right) \\ & = & \phi +  \displaystyle\frac{1}{n} \left(\sum_{i=1}^n \displaystyle\frac{\k_i}{1-\k_i}\right) \psi = \phi +  \left( -1+ \displaystyle\frac{1}{n} \left(\sum_{i=1}^n \displaystyle\frac{1}{1-\k_i}\right)\right) \psi . \end{array}
$$ Now, if we recall the general relation $\lap^g \psi = n{\bf H}$ that holds for any arbitrary spacelike $n$-submanifold of $\L^{n+2}$, we find that $$\phi = \frac{1}{n} \lap^g \psi + \left( 1- \displaystyle\frac{1}{n} \left(\sum_{i=1}^n \displaystyle\frac{1}{1-\k_i}\right)\right) \psi.$$ At last, \eqref{kei} shows that the above equation yields \eqref{lapla}.

Let us compute now $\lap^g \psi$. We shall work at a fixed arbitrary point $x\in U\subset\S^n$. Then there exists an orthonormal basis $\{e_1,\dots,e_n\}$ of $(\S^n,g_0)$ around $x$ such that $$\nabla_{e_i}^{g_0} e_j = 0$$ holds at the point $x$. From now on we will suppress the point $x$ from our notation when possible, as we will always be working at that point.

Let $\{v_0,\dots, v_{n+1}\}$ denote the canonical basis of $\L^{n+2}$, and write $\psi =(\psi_0,\dots, \psi_{n+1})$ in canonical coordinates. Recalling now that $\psi =e^{\rho} (1,x)$, and after expressing the point $x\in \S^n\subset \R^{n+1}\equiv \{v\in \L^{n+2} : v_0 =0\}$ as $x =\sum_{k=1}^{n} x_k v_k$, we compute $$\def\arraystretch{1.8}\begin{array}{lll} \lap^{g_0} \psi & = & \left( \lap^{g_0} (e^{\rho}), \lap^{g_0} (e^{\rho}) x + e^{\rho} \lap^{g_0} x + 2e^{\rho} \displaystyle \sum_{k=1}^n g_0 (\nabla^{g_0} x_k, \nabla^{g_0} \rho) v_k \right)\\ & = & \left( \lap^{g_0} (e^{\rho}), \lap^{g_0} (e^{\rho}) x + e^{\rho} \lap^{g_0} x + 2e^{\rho}  \, \nabla^{g_0} \rho \right) \\ & = & \left( e^{-\rho} \lap^{g_0} (e^{\rho}) \right) \psi + \left(0,  e^{\rho} \lap^{g_0} x + 2e^{\rho}\,   \nabla^{g_0} \rho \right).
\end{array}$$ Now, using that $\lap^{g_0} x = -n x$ and $\lap^{g_0} e^{\rho} = e^{\rho} (\lap^{g_0} \rho + ||\nabla^{g_0} \rho ||_{g_0}^2 )$, we have
 \begin{equation}\label{lagze}
 \lap^{g_0}\psi = (\lap^{g_0} \rho + ||\nabla^{g_0} \rho ||_{g_0}^2 ) \psi + e^{\rho} (0, -nx + 2 \nabla^{g_0} \rho ).
 \end{equation}
Let us introduce now the following notation, for any $Z\in \X(\S^{n})$: $$g_0 (\nabla^{g_0} \psi , Z) := \sum_{k=0}^{n+1} g_0 (\nabla^{g_0} \psi_k,Z) v_k = \sum_{k=0}^{n+1} Z(\psi_k) v_k = Z(\psi)\in \X(\psi)\equiv \X(\phi).$$ With this, we have
\begin{equation}\label{lazei}
\def\arraystretch{1.8}\begin{array}{lll}
g_0(\nabla^{g_0} \psi, \nabla^{g_0} \rho) &=& \left(\nabla^{g_0} \rho \right) (\psi) =\displaystyle \sum_{i=1}^n e_i (\rho) e_i (\psi)  \\ & = &  \displaystyle  \sum_{i=1}^n e_i (\rho) \left(e^{\rho} (e_i (\rho)) (1,x) + e^{\rho} (0,e_i)\right) \\ & = & e^{\rho}  \left( \displaystyle  \sum_{i=1}^n (e_i (\rho))^2 \right) (1,x) + e^{\rho} \left(0, \displaystyle  \sum_{i=1}^n e_i (\rho) e_i\right) \\ & = & ||\nabla^{g_0} \rho ||_{g_0}^2  \, \psi + e^{\rho} \left( 0, \nabla^{g_0} \rho\right).
\end{array}
\end{equation}
Once here, we can recall the usual relation between the Laplacians of two conformal metrics to deduce that, since $g=e^{2\rho} g_0$, we have
 \begin{equation}\label{laconfo}
 \lap^g \psi = e^{-2\rho} \left(\lap^{g_0} \psi + (n-2) g_0 (\nabla^{g_0} \rho,\nabla^{g_0} \psi ) \right).
 \end{equation}
Thus, by \eqref{lagze} and \eqref{lazei}, we obtain from \eqref{laconfo} that
 \begin{equation}\label{lapdef}
 \lap^g \psi = e^{-2 \rho} \left(\lap^{g_0} \rho + (n-1) ||\nabla^{g_0} \rho ||_{g_0}^2 \right)\psi + n e^{-\rho} \left(0, -x + \nabla^{g_0} \rho\right).
 \end{equation}

On the other hand, it is well known that if $g=e^{2\rho} g_0$ is a conformal metric on $\S^n$, the scalar curvature $S(g)$ of $g$ is related to $\rho$ by means of the following elliptic PDE:
 \begin{equation}\label{escnir}
 \Delta^{g_0} \rho + \frac{n-2}{2} ||\nabla^{g_0} \rho||_{g_0}^2 - \frac{n}{2} + \frac{e^{2\rho}}{2(n-1)} S(g) =0.
 \end{equation}
Thus,
 \begin{equation}\label{escnir2}
 e^{-2\rho}\left(\Delta^{g_0} \rho + (n-1) ||\nabla^{g_0} \rho||_{g_0}^2 \right) = -\frac{S(g)}{2(n-1)} + \frac{n e^{-2\rho}}{2} \left(1+||\nabla^{g_0} \rho||_{g_0}^2\right).
 \end{equation}
If we substitute \eqref{escnir2} into \eqref{lapdef} we obtain
 \begin{equation}\label{lapdef2}
 \lap^g \psi = \left(-\frac{S(g)}{2(n-1)} + \frac{n e^{-2\rho}}{2} \left(1+||\nabla^{g_0} \rho||_{g_0}^2\right)\right) \psi + ne^{-\rho} \left(0, -x + \nabla^{g_0} \rho \right).
 \end{equation}
At last, plugging \eqref{lapdef2} into \eqref{lapla} we get \eqref{repfor}, as we wished.
\end{proof}

\begin{remark}\label{catorce}
The parallel flow $\{\phi_t\}_{t\in \R}$ of $\phi$ in $\H^{n+1}\subset \L^{n+2}$ is given by $$\phi_t = \cosh (t) \phi + \sinh (t) (\psi -\phi)= e^{-t} \phi + \sinh (t) \psi.$$ Thus, using \eqref{repfor} we obtain the explicit formula
\begin{equation}\label{remun}
 \phi_{t} = \frac{e^t}{2} \, e^{\rho} \left(1+ \frac{e^{-2\rho}}{e^{2t}} \left( 1+ ||\nabla^{g_0} \rho ||_{g_0}^2 \right)\right) (1,x) + \frac{e^{-\rho}}{e^t} \, (0, -x +\nabla^{g_0} \rho).
 \end{equation}
\end{remark}

\subsection*{A solution for the Christoffel problem in $\H^{n+1}$}

\noindent We are now in the conditions to prove one of our main results.
\begin{teo}\label{main}
Let $\phi:\S^n\flecha \H^{n+1}$ be a solution of the Christoffel problem in $\H^{n+1}$ for the smooth function $C:\S^n\flecha \R_+$. Then its horospherical metric $g_{\8}$ is a solution to the Nirenberg problem in $\S^n$ for the scalar curvature function $S:\S^n\flecha \R$ given by
 \begin{equation}\label{enun1}
 S(x)= n(n-1) (1-2C(x)).
 \end{equation}
Conversely, let $g=e^{2\rho} g_0$ denote a solution to the Nirenberg problem for the scalar curvature function $S:\S^n\flecha \R$. Then there exists $\tau_0 >0$ such that for every $\tau \geq \tau_0$ the map $\phi_{\tau} :\S^n\flecha \H^{n+1}$ given by
 \begin{equation}\label{enun2}
 \phi_{\tau} = \frac{\tau}{2} \, e^{\rho} \left(1+ \frac{e^{-2\rho}}{\tau^2} \left( 1+ ||\nabla^{g_0} \rho ||_{g_0}^2 \right)\right) (1,x) + \frac{e^{-\rho}}{\tau} \, (0, -x +\nabla^{g_0} \rho)
 \end{equation}
is a solution to the Christoffel problem in $\H^{n+1}$ for the smooth function
 \begin{equation}\label{enun3}
C_{\tau} (x)= \frac{1}{2}\left( 1- \frac{S(x)}{\tau^2 n (n-1)}\right).
 \end{equation}
Moreover, the horospherical metric of $\phi_{\tau}$ is actually $g_{\8}= \tau^2 g$.
\end{teo}

\begin{remark}
Theorem \ref{main} proves that the Christoffel problem in $\H^{n+1}$ and the Nirenberg problem on $\S^n$ are \emph{essentially} equivalent problems. On the other hand, we cannot claim that they are completely equivalent. For instance, a scalar curvature function $S(x)$ on $\S^n$ that arises as the horospherical metric of some horospherical ovaloid must satisfy by \eqref{enun1} that $S(x)< n (n-1)$, while this estimate does not hold for general scalar curvature functions on $\S^n$. Alternatively, one can say that given a conformal metric $g=e^{2\rho} g_0$ on $\S^n$, the map $\phi:\S^n\flecha \H^{n+1}$ given by the representation formula \eqref{repfor} will be a solution to the Christoffel problem in $\H^{n+1}$ if and only if it is free of singular points, which is not always the case. We will write down in Corollary \ref{ifanif} an explicit condition for $g$ that is equivalent to the regularity of $\phi$.
\end{remark}

\begin{remark}
Let us explain in more detail the converse in Theorem \ref{main}. First, note that two conformal metrics $g,\tilde{g}$ on $\S^n$ differ just by a dilation (i.e. $\tilde{g} = e^{2t} g$ for some fixed $t\in \R$) if and only if their associated hypersurfaces $\phi,\tilde{\phi}:\S^n\flecha \H^{n+1}$ are parallel (specifically, $\tilde{\phi} =\phi_t$). Thus, given a conformal metric $g= e^{2\rho} g_0$ in $\S^n$, there exists some $\tau_0 >0$ such that if $\tau \geq \tau_0$ the hypersurface $\phi_{\tau}:\S^n\flecha \H^{n+1}$ associated to $g_{\tau} := \tau^2 g$ is everywhere regular, and hence a solution to the Christoffel problem. This is the way that we prevented in Theorem \ref{main} the appearance of singular points for $\phi$.
\end{remark}

\begin{proof}
Let $\phi:\S^n\flecha \H^{n+1}$ denote a solution to the Christoffel problem in $\H^{n+1}$. Then its hyperbolic Gauss map is $G(x)=x$. So, by \eqref{estar} we get that the horospherical metric of $\phi$ is $$g_{\8} = \esiz d\psi,d\psi\esde = e^{2\rho} \esiz dG,dG\esde_{\S^n} = e^{2\rho} g_0,$$ where $g_0$ is the canonical metric of $\S^n$. Hence, the horospherical metric of $\phi$ is globally conformal to $g_0$. Besides, by \eqref{kei} and \eqref{crisdata} we see that \eqref{enun1} holds. This proves the first assertion.

Conversely, let $g= e^{2\rho} g_0$ denote a conformal metric on $\S^n$ with scalar curvature function $S(x)$, and view $\S^n\subset \R^{n+1}$ in the usual way. Consider in addition a positive constant $\tau >0$ and the conformal metric $g_{\tau} := \tau^2 g$, whose scalar curvature function is obviously $S_{\tau} (x)= (1/\tau^2) S(x)$. Then we can construct the map
 \begin{equation}\label{conv1}
 \psi_{\tau} = \tau e^{\rho} (1,x): \S^n\flecha \N_+^{n+1} \subset \L^{n+2} \equiv \L\times \R^{n+1},
 \end{equation}
which has the property that $\esiz d\psi_{\tau},d\psi_{\tau} \esde = \tau^2 e^{2\rho} g_0 =g_{\tau}$. If we let
 \begin{equation*}
 \xi := (0, -x + \nabla^{g_0} \rho),
 \end{equation*}
we may write \eqref{enun2} as
 \begin{equation}\label{conv3}
 \phi_{\tau} = \frac{1}{2} \, \left(1+ \frac{e^{-2\rho}}{\tau^2} \left( 1 + ||\nabla^{g_0} \rho ||_{g_0}^2 \right)\right) \psi_{\tau} + \frac{e^{-\rho}}{\tau} \,\xi.
 \end{equation}
Next, observe that $\esiz \psi_{\tau}, \psi_{\tau}\esde = 0$, and also that $$ \esiz \psi_{\tau},\xi\esde = -\tau e^{\rho} + \tau e^{\rho} \esiz x, \nabla^{g_0} \rho\esde = - \tau e^{\rho},$$ and in the same way $\esiz \xi, \xi \esde = 1+ ||\nabla^{g_0} \rho ||_{g_0}^2.$ Therefore $\esiz \phi_{\tau},\psi_{\tau} \esde =-1$ (thus the first coordinate of $\phi$ in $\L^{n+2}$ is positive), and $$\esiz \phi_{\tau},\phi_{\tau}\esde = -1 - \frac{e^{-2\rho}}{\tau^2} \left( 1+ ||\nabla^{g_0} \rho ||_{g_0}^2 \right) + \frac{e^{-2\rho}}{\tau^2} \left( 1 + ||\nabla^{g_0} \rho ||_{g_0}^2 \right) = -1.$$ This proves that $\phi_{\tau}$ takes its values in $\H^{n+1}\subset \L^{n+2}$. At last, let $Z \in \X (\S^n)$ be a tangent vector field to $\S^n$. Then, using that $Z(x) =Z$, $\esiz x,Z\esde =0$ and $\esiz \nabla^{g_0} \rho, Z\esde = Z(\rho)$, we have by \eqref{enun2} and \eqref{conv1} that $$\def\arraystretch{1.5}\begin{array}{lll} \esiz \phi_{\tau}, Z (\psi_{\tau})\esde &=& \tau \esiz \phi_{\tau}, e^{\rho} Z(\rho) (1,x) + e^{\rho} (0,Z)\esde \\ & = & \esiz (0,-x + \nabla^{g_0} \rho), Z(\rho) (1,x) + (0,Z)\esde \\ & = & - Z(\rho) + \esiz \nabla^{g_0} \rho, Z\esde =0. \end{array}$$ Putting all of this together we have found out that $\phi_{\tau} :\S^n\flecha \H^{n+1}$ is a hypersurface possibly with singular points, but such that at its regular points its associated light cone immersion is $\psi_{\tau}:\S^n\flecha \N_+^{n+1}$. This happens because $\esiz \phi_{\tau},\psi_{\tau}\esde =-1$, $\esiz \psi_{\tau},\psi_{\tau}\esde =0$ and $\esiz d(\phi_{\tau}),\psi_{\tau}\esde =0$, i.e. $\psi_{\tau}$ is normal to $\phi_{\tau}$. In particular, its hyperbolic Gauss map at regular points is $G(x)=x$, and it is horospherically convex at those points. Moreover, the horospherical metric of $\phi_{\tau}$ at regular points is $g_{\8} =\esiz d\psi_{\tau},d\psi_{\tau} \esde = g_{\tau}$, which has the scalar curvature function $S_{\tau}(x) = \tau^{-2} S(x)$. Thus, by \eqref{kei} the mean of the hyperbolic curvature radii of $\phi_{\tau}$ at any regular point is given by $$C_{\tau}= \frac{1}{n} \sum_{i=1}^n \cR_i^{\tau} = \frac{1}{n} \sum_{i=1}^n \frac{1}{1-\k_i^{\tau}} = \frac{1}{2}\left( 1- \frac{S(x)}{\tau^2 n (n-1)}\right),$$ which is exactly \eqref{enun3}. Thus, in order to finish the proof we only need to ensure the existence of some $\tau_0>0$ such that $\psi_{\tau}:\S^n\flecha \H^{n+1}$ is everywhere regular whenever $\tau \geq\tau_0$. First of all, denote $\phi := \phi_1 $ and $\psi :=\psi_1$. Then we may easily observe that
 \begin{equation}\label{flupsi}
\phi_{\tau} = \frac{1}{\tau} \phi + \frac{\tau^2 -1}{2\tau} \psi.
 \end{equation}
Let $U\S^n$ denote the unit tangent bundle of $\S^n$, i.e. $$U\S^n =\{(p,v) \in \S^n\times \S^n : \esiz p,v\esde =0\}.$$ By \eqref{flupsi} we have $d(\phi_{\tau})_p (v)= (1/\tau) d\phi_p (v) + (\tau^2-1)/(2\tau) d\psi_p (v)$ for every $(p,v)\in U\S^n$. Thus, $\phi_{\tau}$ is everywhere regular if and only if $d(\phi_{\tau})_p (v) \neq 0$ for every $(p,v)\in U\S^n$, if and only if
 \begin{equation}\label{flupsi2}
 d\phi_p (v) \neq \frac{1}{2} (1- \tau^2) d\psi_p (v) \hspace{1cm} \text{ for every } (p,v)\in U\S^n.
 \end{equation}
But now, as $U\S^n$ is compact, the sets $$\Lambda := \{ d\psi_p (v) \in \L^{n+2} :  (p,v)\in U\S^n\}, \hspace{1cm} \Omega := \{ d\phi_p (v) \in \L^{n+2}:  (p,v)\in U\S^n\}$$ are compact in $\R^{n+2}\equiv \L^{n+2}$. Moreover, as $d\psi_p (v)\neq 0$ always, we have $0\notin \Lambda$, and thus we may infer the existence of some $r_0 >0$ such that if $r \geq r_0$, then $\Omega\cap r \Lambda = \emptyset$. In particular, there exists some $\tau_0 >0$ such that if $\tau \geq \tau_0$, the condition \eqref{flupsi2} holds. This proves that $\phi_{\tau}$ is everywhere regular if $\tau \geq \tau_0$ and finishes the proof.
\end{proof}

\subsection*{Applications}

As a straightforward consequence of Theorem \ref{main}, we can rephrase into the context of the Christoffel problem in $\H^{n+1}$ all these results on the prescribed scalar curvature problem in $\S^n$. It is our aim now to make explicit some of them.

The following result is a translation into our setting via Theorem \ref{main} of the classical necessary conditions by Kazdan-Warner \cite{KaWa1} and Bourguignon-Ezin \cite{BoEz} on prescribing scalar curvature in $\S^n$.

\begin{cor}[Necessary conditions]
Let $C:\S^n\flecha \R_+$ be a Christoffel function. Then
 \begin{enumerate}
 \item
$C(x) <1/2$ for some $x\in \S^n$.
 \item
If $x_1,\dots, x_{n+1}$ denote the coordinate functions of $\S^n\subset \R^{n+1}$, then $$\int_{\S^n} g_0 \left(\nabla^{g_0} C, \nabla^{g_0} x_i \right) dv_g =0, \hspace{1cm} i=1,\dots , n+1,$$ where $dv_g$ is the volume element of the horospherical metric of the horospherical ovaloid $\phi$. In particular, $C$ cannot be a monotonous function of some coordinate $x_i$.
 \item
More generally, it holds $$\int_{\S^n} X(C) dv_g,$$ for any conformal vector field $X\in \X (\S^n)$.
 \end{enumerate}
\end{cor}

The following very interesting result on the moduli space of solutions to the Nirenberg problem was obtained by Y.Y. Li in \cite{Li1}, as a strong generalization of a previous density result by Bourguignon-Ezin: \emph{smooth scalar curvature functions on $\S^n$ are $C^0$ dense among functions on $\S^n$ that are positive somewhere}. As a consequence of this result and Theorem \ref{main} we have:

\begin{cor}
Let $F:\S^n\flecha \R$ denote a smooth function such that $F(x_0)<1/2$ for some $x_0 \in \S^n$. Then for every $\epsilon >0$ there exists a Christoffel function $C(x):\S^n\flecha \R_+$ and some $t\in \R$ such that $$||F - C_t||_{C^0 (\S^n)} < \epsilon, \hspace{1cm} C_t (x):= \frac{1}{2} - \frac{e^{-2t}}{2} \left(1- 2 C(x)\right).$$ Moreover, the Christoffel function $C(x)$ can be chosen so that $F(x)=C_t (x)$ on the exterior of an arbitrarily small ball of $\S^n$ centered at $x_0$.
\end{cor}

In what respects to sufficient conditions for the Nirenberg problem, we cannot make justice in a few lines to the diversity of results that are known once some technical condition on the scalar curvature function is imposed (a good reference for that is the survey \cite{Li5}). Let us simply say here that all these sufficient conditions can be rephrased one by one in our context to yield sufficient conditions for a function to be a Christoffel function.

We shall nonetheless point out just a couple of sufficient conditions under symmetry assumptions on the Christoffel function. They follow from \cite{Mos,EsSc,ChLi}. For that, let us say that a function $C\in C^{\8} (\S^n)$ is a \emph{generalized Christoffel function} if there exists $t_0\in \R$ such that for every $t\geq t_0$ the map $$C_t (x):= \frac{1}{2} - \frac{e^{-2t}}{2} \left(1- 2 C(x)\right):\S^n\flecha \R$$ is a Christoffel function.

\begin{cor}
Let $C\in C^{\8} (\S^n)$ be a smooth function with $C(x_0)<1/2$ for some $x_0\in \S^n$. Then $C$ is a generalized Christoffel function if it verifies one of the following conditions.
\begin{enumerate}
 \item
$C(x)=C(-x)$ and there is some point $\bar{x}\in \S^n$ such that $C(\bar{x})={\rm min}\, C$ and all derivatives of $C$ up to $n-2$ order vanish at $\bar{x}$.
 \item
$C=C(r)$ is rotationally invariant on $\S^n$, $C'(r)$ changes sign in the region(s) where $C<1/2$, and near any critical point $r_0$ the following \emph{flatness condition} holds: $$C(r)= C(r_0) + a |r-r_0|^{\alfa} + h(|r-r_0|),$$ where $ a\neq 0$, $n-2<\alfa<n$, and $ h'(s)= o (s^{\alfa-1}).$

\end{enumerate}
\end{cor}

\section{Generalized Christoffel problems}

A natural extension of the Christoffel problem in hyperbolic space is to prescribe for a horospherical ovaloid in $\H^{n+1}$ the hyperbolic Gauss map together with a given functional of the hyperbolic curvature radii. Our aim in this section is to show that this problem is equivalent to the question of prescribing a given functional of the eigenvalues of the Schouten tensor for a conformal metric on $\S^n$.

Given a Riemannian metric $g$ on a manifold $M^n$, $n>2$, the \emph{Schouten tensor} of $g$ is the symmetric $(0,2)$-type tensor given by
 \begin{equation}\label{schouten}
 {\rm Sch}_g :=\frac{1}{n-2}\left({\rm Ric}_g - \frac{S(g)}{2(n-1)} \, g \right),
 \end{equation}
where ${\rm Ric}_g$, $S(g)$ stand for the Ricci and scalar curvatures of $g$. Its study is an important issue in conformal geometry, as it represents the non-conformally invariant part of the Riemann curvature tensor. For instance, if $\tilde{g} = e^{2 u} g$ is a conformal metric, then it holds
 \begin{equation}\label{choucon}
 {\rm Sch}_{\tilde{g}} ={\rm Sch}_g - \nabla^{2,g} u + du \otimes du - \frac{1}{2} \, ||\nabla^g u||_g^2 \, g.
 \end{equation}
Once here, let us consider a conformal metric $g= e^{2\rho} g_0$ on $\S^n$, $n>2$. By \eqref{schouten} it follows immediately that ${\rm Sch}_{g_0} = (1/2) g_0$, which is independent of $n$. It then follows by \eqref{choucon} that the Schouten tensor of $g$ is
 \begin{equation}\label{chou2}
 {\rm Sch}_g = - \nabla^{2,g_0} \rho + d\rho \otimes d\rho - \frac{1}{2} \left( -1 + ||\nabla^{g_0} \rho||_{g_0}^2 \right)g_0.
 \end{equation}

The Schouten tensor \eqref{schouten} is not defined for $2$-dimensional metrics. Nevertheless, as \eqref{chou2} makes sense also for $n=2$, we may define naturally the Schouten tensor for conformal metrics on $\S^2$ in a unifying way:

\begin{defi}
Let $g=e^{2\rho} g_0$ denote a conformal metric on $\S^2$. Then its \emph{Schouten tensor} is defined as the symmetric $(0,2)$-type tensor given by \eqref{chou2}.
\end{defi}
Bearing this in mind, we have

\begin{teo}\label{funfors}
Let $\phi:\S^n\flecha \H^{n+1}$ denote a horospherical ovaloid with hyperbolic Gauss map $G(x)=x$, and let $g=e^{2\rho} g_0$ denote its horospherical metric. Then the first and second fundamental forms of $\phi$ at $x\in \S^n$ are given, respectively, by
 \begin{equation}\label{1ffo}
I_{\phi} (e_i,e_j)= \frac{e^{-2\rho}}{4} \left(g(e_i,e_j) - 2\, {\rm Sch}_g (e_i,e_j)\right)^2
\end{equation}
and
 \begin{equation}\label{2ffo}
II_{\phi} = I_{\phi} - \frac{1}{2} \, g + {\rm Sch}_g .
\end{equation}
Here $e_1,\dots, e_n\in T_x \S^n$ is an orthonormal frame with respect to $g_0$, such that $\nabla_{e_i}^{g_0} e_j =0$ at $x$ for every $i,j$.
\end{teo}
\begin{proof}
Let $\{e_1,\dots, e_n\}$ be the above orthonormal basis at $x\in \S^n$, i.e. such that $\nabla_{e_i}^{g_0} e_j =0$ at $x$ for every $i,j$. Then by Theorem \ref{represe} we may write \eqref{repfor} as
 \begin{equation}\label{agrande}
\phi = f (1,x) + e^{-\rho} (0 , \xi) ,
\end{equation}
where $$f:= \frac{e^{\rho}}{2} ( 1+ e^{-2\rho} (1+ ||\nabla^{g_0} \rho ||_{g_0}^2)), \hspace{1cm} \xi := \left(0, -x + \nabla^{g_0} \rho\right) .$$
Then, using that by $\nabla_{e_i}^{g_0} e_j =0$ we have $e_i (e_j)= - \delta_{ij} x$ at $x$, we get
 \begin{equation}\label{bgrande}
e_i (\xi)= \left(0, -e_i + e_i\left( \sum_{j=1}^n e_j(\rho) e_j\right) \right) = \left(0,-e_i - e_i (\rho) x + \sum_{j=1}^n e_i(e_j(\rho)) e_j\right).
\end{equation}
So, by \eqref{agrande} and \eqref{bgrande} we get
 \begin{equation}\label{cgrande}
e_i (\phi) = e_i (f) (1,x) + (f-e^{-\rho}) (0,e_i) - e^{-\rho} \left(0,\sum_{j=1}^n  \{ e_i (\rho) e_j (\rho) - e_i (e_j (\rho))\} e_j \right).
\end{equation}
From \eqref{cgrande} we see directly that $\esiz (1,x),e_i (\phi)\esde =0$. In addition,
 \begin{equation}\label{dgrande}
\esiz e_i (\phi), (0,e_j)\esde = (f-e^{-\rho}) \delta_{ij} - e^{-\rho} \left( e_i (\rho) e_j (\rho) - e_i (e_j (\rho)) \right).
\end{equation}
Thus, by \eqref{cgrande}, \eqref{dgrande} we have
 \begin{equation*}%\label{egrande}
\def\arraystretch{1.8}\begin{array}{lll}
\esiz e_i (\phi), e_j (\phi)\esde &=& e^{-2\rho} (f- e^{-\rho})^2 \delta_{ij} - 2  e^{-2\rho} (f- e^{-\rho}) \left( e_i (\rho) e_j (\rho) - e_i (e_j (\rho)) \right) \\ & & + e^{-2\rho} \left( \displaystyle\sum_{k=1}^n \left(e_i (\rho) e_k (\rho) - e_i (e_k (\rho)) \right)\left(e_j (\rho) e_k (\rho) - e_j (e_k (\rho)) \right)\right).
\end{array}
\end{equation*}
In order to simplify this expression, let $A,B$ denote the $n\times n$ matrices $$ A= (a_{ij}) \hspace{0.3cm} \text{with} \hspace{0.3cm} a_{ij}:= e^{-\rho} (e_i (e_j (\rho)) - e_i (\rho) e_j (\rho)), \hspace{1cm} B:= (f- e^{-\rho}) {\rm Id}_n.$$ Then it holds immediately from the above expressions that $\esiz e_i (\phi), e_j (\phi)\esde$ is the $(i,j)$ entry of the matrix $(A+B)^2$. On the other hand, we may observe that $A$ is the matrix expression of $e^{-\rho} \left(\nabla^{2,g_0} \rho - d\rho \otimes d\rho \right)$ with respect to the basis $\{e_1,\dots, e_n\}$. Thereby, $A+B$ is the matrix expression of $$\frac{e^{\rho}}{2}g_0 + e^{-\rho} \left( \nabla^{2,g_0} \rho - d\rho \otimes d\rho + \frac{1}{2} \left(-1 + ||\nabla^{g_0} \rho ||_{g_0}^2 \right) g_0 \right),$$ which by \eqref{chou2} is nothing but $ e^{-\rho} (g/2 -{\rm Sch}_g ).$ Thus, $$\esiz e_i (\phi),e_j(\phi)\esde =\frac{e^{-2\rho}}{4} \left(g(e_i,e_j) -2\, {\rm Sch}_g (e_i,e_j) \right)^2,$$ what yields \eqref{1ffo}.

In order to compute the second fundamental form, let us first note that $$e_i (\psi)= e^{\rho} e_i (\rho) (1,x) + e^{\rho} (0,e_i),$$ and so by \eqref{dgrande} and \eqref{chou2} we have
 \begin{equation}\label{ggrande}
\esiz d\phi,d\psi\esde = \frac{g}{2} -{\rm Sch}_g.
\end{equation}
Thereby, as $II_{\phi}= -\esiz d\phi,d\eta\esde = \esiz d\phi,d\phi\esde - \esiz d\phi,d\psi\esde$, we obtain \eqref{2ffo} from \eqref{ggrande}.
\end{proof}
As an immediate consequence of this result and Theorem \ref{main}, we have

\begin{cor}\label{ifanif}
A smooth function $C:\S^n\flecha \R_+$ is a Christoffel function if and only if $S(x):= n(n-1)(1-2C(x))$ is the scalar curvature function of some conformal metric $g=e^{2\rho} g_0$ such that $g- 2\, {\rm Sch}_g$ is everywhere positive definite on $\S^n$.
\end{cor}
In \cite{Sch} Schlenker gave a necessary and sufficient condition for a conformally flat metric to be realized as the horospherical metric of some locally horospherically convex hypersurface in $\H^{n+1}$. Thus, Corollary \ref{ifanif} can be seen as an easy consequence of Theorem \ref{main} and \cite{Sch}, although the Schouten tensor ${\rm Sch}_g$ is never mentioned there.

Let $g= e^{2\rho} g_0$ be a conformal metric on $\S^n$, and take $x\in \S^n$. Then we can consider the eigenvectors $v_1,\dots, v_n \in T_x \S^n$ and the eigenvalues $\landa_1,\dots, \landa_n$ of the Schouten tensor ${\rm Sch}_g$ with respect to $g$. Thus, $g(v_i,v_j)=\delta_{ij}$ and ${\rm Sch}_g (v_i,v_j) = \landa_i \delta_{ij}$. With this, Theorem \ref{funfors} yields the following conclusion:

\begin{cor}\label{hipho}
Let $\phi:\S^n\flecha \H^{n+1}$ be a horospherical ovaloid, and let $\{\cR_1,\dots, \cR_n\}$ denote its hyperbolic curvature radii at $x\in \S^n$. Then it holds
 \begin{equation}\label{relac}
 \cR_i = \frac{1}{2} - \landa_i, \hspace{1cm} i=1,\dots, n,
 \end{equation}
where $\{\landa_1,\dots, \landa_n\}$ are the eigenvalues of the Schouten tensor of the horospherical metric $g$ of $\phi$ at $x$. Moreover, the eigendirections of ${\rm Sch}_g$ at $x$ coincide with the principal directions of $\phi$ at $x$.
\end{cor}
\begin{proof}
Let $\{e_1',\dots, e_n'\}$ denote an orthonormal basis of principal directions of $\phi$ at $x$, and define $v_i := 1/(1-\kappa_i) e_i'$. Then $g(v_i,v_j)=\delta_{ij}$ and $$\esiz d\phi (v_i),d\psi (v_j)\esde = \frac{1}{1-\kappa_i} \, \delta_{ij} = \cR_i \, \delta_{ij}.$$
On the other hand, by \eqref{ggrande}, $$\esiz d\phi (v_i),d\psi (v_j)\esde = \frac{1}{2} \delta_{ij} - {\rm Sch_g} (v_i,v_j).$$ Thus ${\rm Sch_g} (v_i,v_j)= (1/2 -\cR_i) \, \delta_{ij}$, which yields \eqref{relac}. The assertion on the eigendirections of ${\rm Sch}_g$ also follows directly.
\end{proof}
Let us point out that from \eqref{relac} and \eqref{crisfun} we get

\begin{equation}\label{sija}
\landa_i = \frac{1}{2} -\cR_i = \frac{1}{2} - \frac{1}{1-\k_i} = -\frac{1+ \k_i}{2(1-\k_i)}.
\end{equation}

\bigskip

\noindent {\bf An interpretation in $\H^{n+1}$ of the Schouten tensor:} Let $\phi:M^n\flecha \H^{n+1}$ denote an oriented hypersurface that is horospherically convex at $p$, and let us assume that $|\k_i (p)|<1$ or $|\k_i (p)|>1$ simultaneously for every $i=1,\dots, n$. Then $\phi$ with its opposite orientation is still horospherically convex at $p$. Thus, its negative hyperbolic Gauss map $G^-$ is a local diffeomorphism at $p$, and we can define the \emph{negative curvature radii} $$\cR_i^* (p) = \frac{1}{|1+ \k_i (p)|}.$$ Once here, we can consider the quotient
 \begin{equation}\label{treses}
 \mathcal{D}_i (p):= \frac{\cR_i (p)}{\cR_i^* (p)} = \frac{|1+ \k_i (p)|}{|1- \k_i (p)|} >0,
  \end{equation}
which measures the variation experimented by the hyperbolic curvature radii of $\phi$ after a change of orientation. We call $\mathcal{D}_i (p)$ the \emph{dilation ratios} of $\phi$ at $p$. This concept does not have an Euclidean counterpart, i.e. it is specific of hyperbolic geometry.

Now, let $\phi:\S^n\flecha \H^{n+1}$ denote a horospherical ovaloid, and assume furthermore that $\k_i <-1$ everywhere. This condition is equivalent (by compactness) to require that $G^-:\S^n\flecha \S^n$ be a global diffeomorphism, and it forces the convexity of $\phi$. We will say that $\phi$ is a \emph{strongly $H$-convex} ovaloid in $\H^{n+1}$.

Let $\{\landa_1,\dots, \landa_n\}$ denote the eigenvalues of the Schouten tensor for the horospherical metric of $\phi$. Then by \eqref{sija} we have
 \begin{equation}\label{unaest}
\landa_i =\frac{1}{2} \mathcal{D}_i.
 \end{equation}
Hence, the eigenvalues $\landa_i$ of the Schouten tensor are identified on a strongly $H$-convex ovaloid as one half of the dilation ratios $\mathcal{D}_i$ of $\phi$.

Alternatively, by \eqref{curra}, we can conclude that
\begin{equation}\label{dosest}
e^{2\varrho_i} = \frac{|1+\k_i|}{|1-\k_i|} = \mathcal{D}_i = 2\landa_i,
\end{equation}
which provides an interpretation of the eigenvalues $\landa_i$ in terms of the classical \emph{contact radii} $\varrho_i= \coth^{-1} (|\k_i|)$.

We will work with a slightly more general notion of contact radii, that is more suitable to the case of non-convex horospherical ovaloids.

\begin{defi}
Let $M^n\subset \H^{n+1}$ denote a horospherical ovaloid with principal curvatures $\k_1,\dots, \k_n$. We define the \emph{signed contact radii} $\delta_1,\dots, \delta_n$ of $M^n$ at $p$ as
 \begin{equation}\label{radel}
 \delta_i (p):= \frac{1+\k_i(p)}{1-\k_i (p)} \hspace{1cm} \text{ {\rm (thus } } |\delta_i (p)| = \mathcal{D}_i (p) \text{{\rm \ if $\k_i(p)\neq -1$)}}.
 \end{equation}
\end{defi}

\vspace{0.2cm}

\noindent {\bf The generalized Christoffel problem:}

\begin{quote}
Let $\Gamma :\S^n\flecha \R$ be a smooth function. Find out if there exists a hypersurface $\phi:\S^n\flecha \H^{n+1}$ (necessarily a horospherical ovaloid) whose hyperbolic Gauss map is $G(x)=x$ for every $x\in \S^n$, and such that
 \begin{equation}\label{gcris}
\mathcal{F}(\cR_1, \dots, \cR_n) = \Gamma ,
 \end{equation}
where here $\mathcal{F}(\cR_1, \dots, \cR_n)$ is a prescribed functional of the hyperbolic curvature radii $\cR_i$ of the hypersurface $\phi$.
\end{quote}

\begin{remark}
By \eqref{crisfun} it is obvious that prescribing a relation of the form \eqref{gcris} for the hyperbolic curvature radii $\cR_i >0$ can also be seen as prescribing a relation of the form $\mathcal{W} (\k_1,\dots, \k_n)=\Gamma$ for the principal curvatures $\k_1,\dots, \k_n$, with the restriction $\k_i \in (-\8,1)$. The novelty of the generalized Christoffel problem with respect to other previous works on hypersurfaces in space forms satisfying $\mathcal{W} (\k_1,\dots, \k_n)=\Gamma$ is that we are also prescribing here the hyperbolic Gauss map $G$.
\end{remark}

The problem in this generality is quite hopeless to be solved. Nonetheless, we can conclude from Corollaries \ref{ifanif} and \ref{hipho} the following:

\begin{cor}\label{gcr}
Let $\phi:\S^n\flecha \H^{n+1}$ denote a solution to the generalized Christoffel problem, let $g$ denote its horospherical metric, and define
 \begin{equation}\label{fg}
\mathcal{G} (t_1,\dots, t_n):= \mathcal{F} (1/2- t_1,\dots, 1/2- t_n).
 \end{equation}
Then the eigenvalues $\landa_i$ of ${\rm Sch}_g$ verify $\landa_i <1/2$ and $\mathcal{G} (\landa_1,\dots, \landa_n)=\Gamma (x).$

Conversely, if $g=e^{2\rho} g_0$ is a conformal metric on $\S^n$ such that the eigenvalues $\landa_i$ of its Schouten tensor verify $\landa_i <1/2$ and $\mathcal{G} (\landa_1,\dots, \landa_n)=\Gamma (x)$ for a certain functional $\mathcal{G}$, then there exists a solution to the generalized Christoffel problem for $\mathcal{F}$ given by \eqref{fg}, and whose horospherical metric is $g$.
\end{cor}

This result steps the path for studying the generalized Christoffel problem in special cases, as the behavior of the eigenvalues of the Schouten tensor is a more studied subject. We discuss next some interesting points of the correspondence given by Corollary \ref{gcr}.

\subsection*{The $k$-th symmetric problems on $\S^n$ and $\H^{n+1}$}

The problem of prescribing a certain functional of the eigenvalues of the Schouten tensor for a conformal metric on a Riemannian manifold has received several contributions in the last few years, see \cite{Cha,CGY,CGY2,Via,JLX} and specially the book \cite{Gua} and references therein. In particular, special attention has been paid to the study of the $\sigma_k$\emph{-curvatures} of a Riemannian metric $g$, that can be defined as follows.

First, let $S_k (x_1, \dots, x_n)$ denote the $k$-th elementary symmetric function of $(x_1,\dots, x_n)$, i.e.
 \begin{equation}\label{ksim}
  S_k (x_1, \dots, x_n) = \sum x_{i_1} \cdots x_{i_k},
 \end{equation}
where the sum is taken over all strictly increasing sequences $i_1, \dots, i_k$ of the indices from the set $\{1,\dots , n\}$. Then, given a metric $g$ on a manifold $M^n$, the $\sigma_k$\emph{-curvature} of $g$ is defined as $S_k (\landa_1,\dots, \landa_n)$, where $\landa_i$ are the eigenvalues of the Schouten tensor of $g$.

The problem of finding a conformal metric on $\S^n$ with a prescribed $\sigma_k$-curvature has been specially treated, as a generalization of the Nirenberg problem. By Corollary \ref{gcr} and \eqref{unaest}, this $\sigma_k$-Nirenberg problem is equivalent (up to dilations) to a natural problem for horospherical ovaloids in $\H^{n+1}$. Namely, the generalized Christoffel problem in which a $k$-th symmetric elementary function of the signed dilation ratios $\delta_i$ in \eqref{radel} is prescribed. Consequently, the local estimates in \cite{Gua} or the Kazdan-Warner type obstructions in \cite{Han} can be translated into analogous results for the hyperbolic setting. We omit the details, as the process is clear from Corollary \ref{hipho}.

The $\sigma_n$-problem is of special interest in $\H^{n+1}$. Recall that the contact radii $\varrho_i$ in \eqref{curra} can be seen to some extent as analogous quantities to the classical curvature radii in $\R^{n+1}$. The negative point was that they only make sense for strongly $H$-convex ovaloids in $\H^{n+1}$, which is quite restrictive. Nonetheless, by \eqref{dosest} we can see that
 \begin{equation}\label{conw}
\Lambda := \frac{1}{n} \sum_{i=1}^n \varrho_i =\frac{\log (2^n \sigma_n)}{2 n}.
 \end{equation}
Hence, by Corollary \ref{gcr}, the $\sigma_n$-Nirenberg problem on $\S^n$ is equivalent (up to a regularity condition) to the \emph{contact} Christoffel problem for strongly $H$-convex ovaloids in $\H^{n+1}$, in which the hyperbolic Gauss map together with the mean $\Lambda$ of the contact radii $\varrho_i$ are prescribed.

On the other hand, it is a very natural question to consider {\bf the $k$-th symmetric Christoffel problem in $\H^{n+1}$,} i.e. the generalized Christoffel problem in $\H^{n+1}$ with the specific choice $\mathcal{F} =S_k$.

It is clear that if we take $k=1$ we get the Christoffel problem in $\H^{n+1}$. The $k$-th symmetric Christoffel problem is a natural analogue in $\H^{n+1}$ of the \emph{Christoffel-Minkowski problem} for ovaloids in $\R^{n+1}$, which prescribes the Gauss map as well as a $k$-th symmetric function of the principal curvature radii. Nonetheless, in this hyperbolic setting, the Minkowski problem (i.e. prescribing the hyperbolic Gauss map along with the Gauss-Kronecker curvature) does not appear as any of the $k$-th symmetric problems.

Now, by \eqref{relac} we have $$S_k (\cR_1,\dots, \cR_n)= \sum_{j=0}^k \left(\def\arraystretch{0.5}\begin{array}{c} n-j \\ \\ k -j
\end{array}\right)(-1)^j  \, 2^{j-k} S_j (\landa_1,\dots, \landa_n).$$

So, by this relation and Corollary \ref{gcr}, the above $k$-th symmetric problem in $\H^{n+1}$ is equivalent (up to a regularity condition) to the Nirenberg-type problem in $\S^n$ in which the following linear combination of the $\sigma_k$-curvatures is prescribed: $$\sum_{j=0}^k c_j \, \sigma_j (g) =\Gamma, \hspace{1cm} c_j := \left(\def\arraystretch{0.5}\begin{array}{c} n-j \\ \\ k -j
\end{array}\right)(-1)^j \, 2^{j-k}.$$ Thus, a natural situation in $\H^{n+1}$ motivates the (seemingly unexplored) problem of prescribing a linear combination of the $\sigma_k$-curvatures for conformal metrics on $\S^n$.

%usar que $$S_k (I+D) =\sum_{j=0}^k \left(\def\arraystretch{0.5}\begin{array}{c} n-j \\ k -j
%\end{array}\right) \, S_j (D), \hspace{0.5cm} S_k (c D)= c^k S_k (D), \hspace{0.4cm} S_k (D^{-1})= \frac{S_{n-k} (D)}{S_n (D)}.$$

Let us finally point out that there are many other interesting problems for hypersurfaces in $\H^{n+1}$ that can be formulated in terms of geometric PDEs problems by means of Corollary \ref{hipho} and Corollary \ref{gcr}. For instance, Dirichlet problems at infinity for horospherically convex hypersurfaces with prescribed Weingarten curvatures, or Christoffel-Minkowski type problems for horospherically convex hypersurfaces with $1$ or $2$ points at the ideal boundary of $\H^{n+1}$. So, an important point of the construction developed in this section is that it motivates new interesting problems for conformally invariant PDEs in terms of the Schouten tensor.

\section{Weingarten hypersurfaces in $\H^{n+1}$}

An oriented hypersurface in a model space $\R^{n+1}$, $\S^{n+1}$, $\H^{n+1}$ is a \emph{Weingarten hypersurface} if there exists a non-trivial relation
 \begin{equation}\label{w}
 \cW (\k_1,\dots, \k_n)=c\in \R
 \end{equation}
between its principal curvatures $\k_1,\dots, \k_n$. In particular, totally umbilical hypersurfaces are Weingarten hypersurfaces. An important problem in hypersurface theory is to establish for which functionals $\cW$ the only examples that satisfy \eqref{w} together with a global condition (usually compactness) are totally umbilical.

Our next objective is to exhibit a wide family of functionals $\cW$ for which the only compact hypersurfaces verifying \eqref{w} are totally umbilical round spheres. This result will be obtained as a consequence of Corollary \ref{hipho} and a deep result by A. Li and Y.Y. Li \cite[Corollary 1.6]{LiLi1} on the Schouten tensor of conformal metrics on $\S^n$.

There are two remarkable facts regarding our result. To start, it seems to be the
first example of a wide family of Weingarten functionals for which round spheres can
be characterized as the only compact examples satisfying \eqref{w}. Previous results
in this line for specific Weingarten functionals in $\H^{n+1}$ may be found, for
instance, in \cite{Ale,IB,CD,SEK,HR3} and references therein. The other remarkable
feature of our result is that we are not assuming \emph{a priori} that the
hypersurface is embedded. This additional topological hypothesis appears in most of
the previous geometric works on the characterization of compact Weingarten
hypersurfaces in model spaces.

In order to state the result, some notation is needed. Let us first of all denote $\Omega := \{(x_1,\dots, x_n)\in \R^{n} : x_i <1\}$, and consider the involutive diffeomorphism $ \cT:\Omega \flecha \Omega$ given by
 \begin{equation}\label{w1}
  \cT (x_1,\dots, x_n)= \left(\frac{x_1 +1}{x_1-1}, \dots, \frac{x_n+1}{x_n-1}\right).
   \end{equation}
Consider in addition
 \begin{equation}\label{w2}
 \Gamma \subset \R^n \text{ an open convex symmetric cone with vertex at the origin}
 \end{equation}
such that
\begin{equation}\label{w3}
\Gamma_n \subset \Gamma \subset \Gamma_1,
 \end{equation}
 where $\Gamma_n := \{(x_1,\dots, x_n ) : x_i >0\}$,  and $\Gamma_1 := \{(x_1,\dots, x_n): \sum_{i=1}^n x_i >0 \}.$ Here \emph{symmetric} means symmetric with respect to the variables $(x_1,\dots, x_n)$.

Let us define next
 \begin{equation}\label{w4}
 \Gamma^* :=\cT(\Gamma \cap \Omega)\subset \Omega.
 \end{equation}
It is clear that $\Gamma_n^* \subset \Gamma^* \subset \Gamma_1^*$, where $$\Gamma_n^* := \{(x_1,\dots, x_n): x_i <-1\}, \hspace{0.5cm} \Gamma_1^* := \left\{(x_1,\dots, x_n): \frac{1}{n}\sum_{i=1}^n \frac{1}{1-x_i} < \frac{1}{2} \right\}.$$ Consider now a real function $\cW (x_1,\dots, x_n)$ satisfying the following conditions:
 \begin{equation}\label{w5}
 \cW \in C^1 (\Gamma^*)\cap C^0 \left(\overline{\Gamma^*}\right) \text{ is symmetric with respect to $x_1,\dots, x_n$.}
 \end{equation}

 \begin{equation}\label{w6}
 \cW = 0 \text{ on } \parc \Gamma^* \text{ and } \parc \cW / \parc x_i <0 \text{ in } \Gamma^*.
 \end{equation}

\begin{equation}\label{w7}
 s\cW (x) = \cW(\cT (s (\cT(x)))) \text{ for every } x\in \Gamma^* .
 \end{equation}
In this last condition $s>0$ is arbitrary except for the fact that the right hand side of \eqref{w7} has to be well defined. Let us observe that if we denote
 \begin{equation}\label{w8}
 f(x):= \cW (\cT (2x)) :\Gamma \cap \{(x_1,\dots, x_n) : x_i < 1/2\} \flecha \R,
 \end{equation}
then \eqref{w7} simply means that $f(s x)= s f(x)$. This fact allows us to extend $f$ by homogeneity to the whole cone $\Gamma$. With this, such extension (still denoted by $f$) verifies:
 \begin{equation}\label{w9}
 f\in C^1 (\Gamma)\cap C^0 \left(\overline{\Gamma}\right) \text{ is symmetric with respect to $x_1,\dots, x_n$.}
 \end{equation}

 \begin{equation}\label{w10}
 f = 0 \text{ on } \parc \Gamma \text{ and } \parc f / \parc x_i >0 \text{ in } \Gamma.
 \end{equation}

\begin{equation}\label{w11}
 f(sx)= s f(x) \text{ for every } x\in \Gamma , s>0.
 \end{equation}

Under the conditions \eqref{w2}, \eqref{w3},  \eqref{w9}, \eqref{w10}, \eqref{w11} for $(\Gamma,f)$, A. Li and Y.Y. Li proved in \cite{LiLi1} that if $g= e^{2\rho} g_0$ is a conformal metric on $\S^n$ whose eigenvalues of its Schouten tensor verify
 \begin{equation}\label{w12}
 f(\landa_1,\dots, \landa_n) =1, \hspace{0.5cm} (\landa_1,\dots, \landa_n )\in \Gamma,
 \end{equation}
then $g$ differs from $g_0$ by, at most, a dilation and a conformal isometry of $\S^n$.

As a consequence of this result in \cite{LiLi1} and of our previous discussion, we can conclude the following:
 \begin{teo}\label{wt1}
Let us define $(\Gamma^*,\cW)$ by the conditions \eqref{w2} to \eqref{w7}, and let $M^n\subset \H^{n+1}$, $n>2$, denote an immersed oriented compact hypersurface in $\H^{n+1}$ such that
\begin{equation}\label{w13}
\cW (\k_1,\dots, \k_n)=1
\end{equation}
holds for its principal curvatures $\k_i$. Then $M^n$ is a totally umbilical round sphere.
 \end{teo}
\begin{proof}
Since $\Gamma^* \subset \Omega$, we infer that $\k_i <1$ always holds, and thus $M^n$ is a canonically oriented horospherical ovaloid in $\H^{n+1}$. Thereby, $M^n$ is diffeomorphic to $\S^n$ and its horospherical metric $g_{\8}$ is conformal to the canonical metric of $\S^n$. Moreover, by \eqref{crisfun} and \eqref{relac} we can express the eigenvalues of the Schouten tensor of $g_{\8}$ as
 \begin{equation}\label{w14}
 2\landa_i = \frac{\k_i +1}{\k_i -1}.
 \end{equation}
Let us recall at this point that the conditions imposed to $(\Gamma^*, \cW)$ ensure that $(f,\Gamma)$ given by \eqref{w2}, \eqref{w3}, \eqref{w4}, \eqref{w8} are in the conditions of \cite[Corollary 1.6]{LiLi1}. Besides, we get from \eqref{w7}, \eqref{w13}, \eqref{w14} that $$f(\landa_1,\dots, \landa_n)=1, \hspace{0.5cm} (\landa_1,\dots, \landa_n)\in \Gamma.$$ Consequently, $g_{\8}$ is a metric of constant (positive) curvature, and all eigenvalues $\landa_i$ are constant and equal. By \eqref{w14} we conclude that $M^n$ is a totally umbilical round sphere.
 \end{proof}

\noindent {\bf An important remark:} observe that if the Weingarten functional $\cW (x_1,\dots, x_n)$, which is defined on $\Gamma^*$, admits an extension $\widetilde{\cW}$ to some larger domain $\widetilde{\Gamma}\subset \R^n$, then the condition \eqref{w13} on an immersed oriented hypersurface $M^n\subset \H^{n+1}$ simply means that $$\widetilde{\cW} (\k_1,\dots, \k_n)=1$$ and there exists some $p\in M^n$ with $(\k_1(p),\dots, \k_n(p))\in \Gamma^*$. This happens because $\cW$ vanishes identically on $\parc \Gamma^*$.
So, if $M^n$ is compact, the condition $(\k_1(p),\dots, \k_n(p))\in \Gamma^*$ is not restrictive because, by a change of orientation if necessary, there is some $p\in M^n$ with $\k_i (p)<-1$ for every $i$. That is, $(\k_1(p),\dots, \k_n(p))\in \Gamma_n^*\subset \Gamma^*$.

\vspace{0.7cm}

Let us now consider a particular case of Theorem \ref{wt1} that is specially important. Let $S_k (x_1,\dots, x_n)$ denote the elementary $k$-th symmetric function \eqref{ksim}, and let $\Gamma_k\subset \R^n$ denote the connected component of  $$\{(x_1,\dots, x_n)\in \R^n: S_k (x_1, \dots, x_n) >0\}$$ that contains the positive cone $\Gamma_n$. It is then known (see \cite{LiLi1,LiLi2,Li3}) that $(\Gamma_k,f_k)$ verify conditions \eqref{w2}, \eqref{w3},  \eqref{w9}, \eqref{w10}, \eqref{w11}, where $$f_k (x_1,\dots, x_n):= \left( S_k (x_1,\dots, x_n)\right)^{1/k}.$$ Thereby, we obtain from Theorem \ref{wt1} and \eqref{w14}

\begin{cor}\label{lwcor}
Let $M^n\subset \H^{n+1}$, $n>2$, denote a horospherical ovaloid, and assume that a $k$-th symmetric elementary function of its signed contact radii is constant, i.e.
 \begin{equation}\label{lwrar}
 S_k (\delta_1,\dots, \delta_n) = c>0, \hspace{1cm} \delta_i := \frac{\k_i+1}{\k_i-1}.
 \end{equation}
Then $M^n$ is a totally umbilical round sphere.
\end{cor}
The particular case $k=n$ tells (by \eqref{conw}) that \emph{any strongly $H$-convex ovaloid in $\H^{n+1}$ with constant mean of their contact radii, i.e. $\sum_{i=1}^n \varrho_i = c>0$, is a round sphere}.

The above Corollary is relevant in the context of linear Weingarten hypersurfaces. Let us recall that for an oriented hypersurface $M^n\subset \H^{n+1}$ with principal curvatures $\k_1,\dots, \k_n$, the $r$-th mean curvature function is defined as $$\left(\def\arraystretch{0.5}\displaystyle \begin{array}{c} n \\ \\ r
\end{array}\right) H_r := S_r (\k_1,\dots, \k_n), \hspace{1cm} 1\leq r\leq n.$$ In particular, $H=H_1$ is the mean curvature and $K=H_n$ is the Gauss-Kronecker curvature. With these notations, we recall the following
 \begin{defi}
An immersed oriented hypersurface $M^n\subset \H^{n+1}$ is a \emph{linear Weingarten hypersurface} if there are constants $c_0,\dots, c_n\in \R$ such that
 \begin{equation}\label{lwf}
 \sum_{i=1}^n c_r H_r = c_0.
 \end{equation}
 \end{defi}

Now, observe that, for any oriented horospherically convex hypersurface, the equality \eqref{lwrar} can be rewritten as (we use here $r$ instead of $k$ for clarity) $$\sum_{i_1<\cdots <i_r} \left(\frac{\k_{i_1} + 1}{\k_{i_1}-1}\right) \cdots \left(\frac{\k_{i_r} + 1}{\k_{i_r}-1}\right) = c  ,$$ or equivalently, $$\sum_{i_1<\cdots <i_r} (\k_{i_1} + 1)\cdots (\k_{i_{r}}+1) (\k_{j_1} - 1)\cdots (\k_{j_{n-r}}-1) = c \,  (\k_1-1)\cdots (\k_n -1),$$ where we have labelled $\{i_1,\dots, i_r,j_1,\dots, j_{n-r}\}= \{1,\dots, n\}$. This proves that the Weingarten relation given by \eqref{lwrar} can be actually written as a linear Weingarten relation
$$ \sum_{i=1}^n c_r H_r = c_0$$ for adequate constants $c_0,\dots, c_n$ depending on $c$ and $r$.

Thus, as a consequence of Corollary \ref{lwcor} we can conclude that \emph{the only
compact oriented immersed hypersurfaces in $\H^{n+1}$ that satisfy \eqref{lwf} on
$\Gamma_r^*$ (given by \eqref{w4} in terms of $\Gamma_r$) for the above adequate
constants $c_j$ are totally umbilical round spheres}. This last conclusion was
obtained in \cite{FeRo} for $k=1$.

\begin{remark}
A problem of interest in hypersurface theory is to classify the elements $(c_0,\dots, c_n)\in \R^{n+1}$ for which the only compact linear Weingarten hypersurfaces verifying \eqref{lwf} are round spheres. For that, one obviously needs some condition on the $c_i$'s in order to ensure that round spheres belong to this class. Again, this problem is tightly linked by Corollary \ref{gcr} to the one of determining when a conformal metric $g=e^{2\rho} g_0$ on $\S^n$ whose $\sigma_k$-curvatures verify a relation of the type $$\sum_{k=1}^n c_k \sigma_k = c_0,\hspace{1cm} (c_0,\dots, c_n)\in \R^{n+1},$$ differs from $g_0$ at most by a dilation and a conformal isometry of $\S^n$.
\end{remark}

\subsection*{Weingarten hypersurfaces with one end}

It is possible to extend partially the above arguments to the case of non-compact Weingarten hypersurfaces. Again, the key is the existence of some highly developed theorems by Y.Y. Li \cite{Li3,Li4} on the characterization of solutions globally defined on $\R^n$ of some conformally invariant PDEs involving the Schouten tensor. Such results by Li generalize other previous works, like \cite{CGS,CGY2,LiLi1,LiLi2}. For that we introduce the following notions.

\begin{defi}
A manifold $M^n$ is said to have \emph{one end} if it is homeomorphic to a compact manifold $K^n$ with one point removed, i.e. $M^n \simeq K^n \setminus \{p\}.$
\end{defi}

\begin{defi}
Let $M^n\subset \H^{n+1}$ denote an immersed oriented hypersurface with one end. We say that this end is \emph{regular} (or that $M^n$ has one \emph{regular} end) provided the hyperbolic Gauss map $G$ of $M^n\equiv K^n \setminus \{p\}$ extends continuously to $K^n$.
\end{defi}
The most basic examples of non-compact hypersurfaces in $\H^{n+1}$ with one regular
end are horospheres (with the outer orientation). The following result shows that,
for a large class of Weingarten functional equations $\cW (\k_1,\dots, \k_n)=0$, the
only Weingarten hypersurfaces with one regular end are horospheres. We are not
assuming here completeness or embeddedness of the hypersurface. Previous
characterizing results in $\H^{n+1}$ for horospheres can be found in \cite{FeRo,Cur}.

\begin{teo}
Let $(\Gamma^*, \cW)$ be given by \eqref{w2} to \eqref{w6} and such that $\cW >0$ on $\Gamma^*$. Let $M^n\subset \H^{n+1}$, $n>2$, denote an immersed oriented hypersurface with one regular end, and whose principal curvatures $\k_1,\dots, \k_n$ verify
 \begin{equation}\label{whor}
 (\k_1,\dots, \k_n)\in \parc \Gamma^*,\hspace{1cm} \text{ {\rm and so }} \cW(\k_1,\dots, \k_n)=0.
 \end{equation}
Then $M^n$ is a horosphere (with its outer orientation).
\end{teo}
\begin{proof}
Let us start by observing that if $G:M^n\flecha \S^n$ is the hyperbolic Gauss map of
$M^n$, then $G$ is a local diffeomorphism at every point of $M^n$ (since from
\eqref{whor} we know that $\k_i <1$). In addition, as the end of $M^n \equiv
K^n\setminus \{p\}$ is regular we can assure the existence of a continuous map
$\bar{G}:K^n\flecha \S^n$ with $\bar{G}|_{M^n} = G$.

Let us denote $x_0 := \bar{G} (p)\in \S^n$, and $N:= G^{-1} (x_0) \subset M^n$. Obviously, $N$ consists only of isolated points, since $G$ is a local diffeomorphism. Thus $M^n\setminus N$ is a smooth manifold, and $$G:M^n \setminus N\flecha \S^n \setminus \{x_0\}$$ is a covering map. Therefore $M^n\setminus N$ is simply connected, and a simple topological argument immediately implies that $N$ has to be empty. Then, we can infer that $$G:M^n\flecha \S^n\setminus \{x_0\}$$ is a global diffeomorphism. In particular $M^n$ is homeomorphic to $\R^n$, and so we can see it as an immersion $\phi :\R^n\flecha \H^{n+1}$ with hyperbolic Gauss map $G= \pi^{-1} (x)$, where $\pi:\S^n \setminus \{\text{ north pole }\} \flecha \R^n$ is the stereographic projection.

Once here, the rest of the proof is basically the same as the one of Theorem \ref{wt1}, but this time using the main results in the works \cite{Li3,Li4}.
\end{proof}

\subsection*{Results on the Schouten tensor}

We use now the interrelation between hypersurfaces in $\H^{n+1}$ and conformal
metrics on $\S^n$ in the opposite direction, i.e. we will induce results from
hypersurface theory into conformal geometry results.

First, we will prove a general duality for locally conformally flat metrics on a manifold
$M^n$. It is an example of how a geometrically simple transformation (a change from
positive to negative hyperbolic Gauss map in our situation) can yield a non-trivial
superposition formula for a geometric PDE.

\begin{teo}\label{inverse}
Let $(M^n,g)$, $n>2$, denote a locally conformally flat Riemannian manifold, and let $\landa_1,\dots, \landa_n \in C^{\8} (M^n)$ be the eigenvalues of its Schouten tensor.

If none of the $\landa_i$'s vanish on $M^n$, then there exists a new locally conformally flat Riemannian metric $g^*$ on $M^n$ whose Schouten tensor eigenvalues $\landa_i^*\in C^{\8} (M^n)$ are $$\landa_i^* = \frac{1}{\landa_i}.$$

%, \ \dots, \ \landa_n^* = \frac{1}{\landa_n}.$$
\end{teo}
\begin{proof}
Take $p_0\in M^n$ an arbitrary point. Then there exists some relatively compact neighbourhood $U\subset M^n$ of $p_0$ and a conformal embedding $f:U\flecha \S^n$. Moreover, by \cite{Kui}, the map $f$ on $U$ is unique up to conformal transformations of $\S^n$.

Consider now the dilated metric $g_t := e^{2 (\rho +t)} g_0$ on $U$, where here $t>0$ is a positive constant. As this is just a dilation of $g$ of factor $e^{2t}$, the eigenvalues of the Schouten tensor of $g_t$ are $\landa_i^t =e^{-2t} \landa_i \neq 0$. Thereby, as $U$ is relatively compact, there exists $t_0 >0$ such that $\landa_i^t \neq 1/2$ everywhere on $U$ if $t\geq t_0$.

Consider once here the map $$\psi_t (x):= e^{\rho (x)+t} \, (1,x):U\subset \S^n\flecha \N_+^{n+1},$$ which is an immersion with induced metric $g_t$. It follows then from the proof of Theorem \ref{main} that the map $\phi_{t} :U\subset\S^n\flecha \H^{n+1}$ given by
 \begin{equation}\label{inv1}
 \phi_{t} = \frac{e^{\rho + t}}{2} \left(1+ e^{-2(\rho+t)} \left( 1+ ||\nabla^{g_0} \rho ||_{g_0}^2 \right)\right) (1,x) + e^{-(\rho +t)} \, (0, -x +\nabla^{g_0} \rho)
 \end{equation}
has $\psi_t$ as its associated light cone immersion for a specific orientation, i.e. $\esiz d\phi_t,d\psi_t\esde=0$ and $\esiz \phi_t,\psi_t\esde = -1$. Moreover, by \eqref{1ffo}, the map $\phi_t$ is everywhere regular, and hence an immersed hypersurface in $\H^{n+1}$. We shall orient $\phi_t$ so that its unit normal $\eta_t$ is the one that verifies $\psi_t = \phi_t + \eta_t$.

Now, by \eqref{curprin} and \eqref{sija} we find that the principal curvatures $\k_i^t$ of $\phi_t$ verify
 \begin{equation}\label{inv2}
 e^{-2t} \landa_i = \frac{1}{2} - \frac{1}{1-\k_i^t},
 \end{equation}
and thus $\k_i^t \neq \pm 1$ at every point in $U$ (recall that, by hypothesis, the eigenvalues $\landa_i$ never vanish).

This implies that the negative Gauss map of $\phi_t$ is a local diffeomorphism, since by Lemma \ref{hc} and after a change of orientation, this condition for $$G_t^- :=[\phi_t - \eta_t]:U\subset \S^n\flecha \S^n$$ is equivalent to $\k_i^t \neq -1$ for every $i\in \{1,\dots, n\}$. Moreover,
 \begin{equation}\label{inve3}
g_t^* := \esiz d(\phi_t - \eta_t),d(\phi_t - \eta_t)\esde
 \end{equation}
is a new regular metric on $U$, that is conformally flat since it is the induced metric of an immersion in the light cone. In other words, $g_t^*$ is just the horospherical metric of $\phi_t$ after a change of orientation, which as we already know, changes completely the geometry of the hyperbolic Gauss map.

Now, by \eqref{sija}, the eigenvalues $\landa_i^{t,*}$ of the Schouten tensor of $g_t^*$ verify the relation $$\landa_i^{t,*} = \frac{1}{2} - \frac{1}{1-\k_i^t}.$$ From this expression and \eqref{inv2} we get $$\landa_i^{t,*} = \frac{e^{2t}}{4\landa_i}.$$ It follows hence directly from here  that the metric
 \begin{equation}\label{inv3}
 g^*:= e^{2t} g_t^*
 \end{equation}
defined on $U$ has the Schouten tensor eigenvalues $\landa_i^* = 1/(4\landa_i)$.

Let us check finally that the definition of the conformally flat metric $g^*$ on $U$ does not depend on the choices of the positive constant $t>0$ and of the conformal embedding $f$.

Let $t_1,t_2 \geq t_0 >0$ be two positive real numbers, and denote by $$\psi_{t_1}(x) =e^{\rho(x) +t_1} \, (1,x), \hspace{1cm} \psi_{t_2} (x)= e^{\rho (x)+ t_2} \, (1,x)$$ the light cone immersions associated to them. Then by Remark \ref{catorce} we have $$\phi_{t_2} = e^{-(t_2-t_1)} \phi_{t_1} + \sinh (t_2-t_1) \psi_{t_1}.$$ Thus, $\phi_{t_2} - \eta_{t_2} = e^{t_1-t_2} (\phi_{t_1} -\eta_{t_1}).$ But this equality together with \eqref{inve3} show that $e^{2t_2} g_{t_2}^* = e^{2t_1} g_{t_1}^*.$ This indicates that the definition of the metric $g^*$ in \eqref{inv3} is independent of the value of the constant $t>0$.

On the other hand, as we said, the conformal embedding $f$ of $U$ into $\S^n$ is unique up to conformal transformations of $\S^n$. This implies that $\psi_t, \phi_t$ and thus $\phi_t - \eta_t$ are unique up to isometries of $\N_+^{n+1}$ or $\H^{n+1}$. In particular, $g^*$ is independent of the choice of the conformal embedding $f:U\flecha \S^n$.

These independence properties of the metric $g^*$ clearly imply that if $U_1,U_2$ are two relatively compact neighborhoods of $M^n$, and if $g_1^*, g_2^*$ denote there respective associated metrics defined via \eqref{inv3}, then $g_1^* = g_2^*$ on $U_1\cap U_2$. In this way, the metric $g^*$ is a locally conformally flat Riemannian metric globally defined on $M^n$. Moreover, its Schouten tensor eigenvalues at an arbitrary point are $1/(4\landa_i)$. Finally, a dilation of this metric $g^*$ by an adequate constant gives the desired metric.
\end{proof}
Let us point out regarding this result that the conformally flat metrics $g,g^*$ are not in general conformal to each other on $M^n$.

The final result of this work is a two-dimensional analogue of the results in
\cite{LiLi1} regarding conformal metrics on $\S^n$. It is remarkable that in this
$2$-dimensional case the hypothesis are much weaker.

\begin{teo}
Let $g= e^{2\rho} g_0$ denote a conformal metric on $\S^2$, and let
$\landa_1,\landa_2$ be the eigenvalues of its $2$-dimensional Schouten tensor, i.e.
of $${\rm Sch}_g := - \nabla^{2,g_0} \rho + d\rho \otimes d\rho - \frac{1}{2} \left(
-1 + ||\nabla^{g_0} \rho||_{g_0}^2 \right)g_0.$$ Assume that
$\cW(\landa_1,\landa_2)=0$ for a smooth function $\cW:D\subset
\R^2\flecha\R$ such that
$$\frac{\parc\cW}{\parc \landa_1}\frac{\parc \cW}{\parc \landa_2} >0 \hspace{1cm}
\text{ whenever } \landa_1=\landa_2.$$ Then $\landa_1\equiv\landa_2$ and $g=\Phi^*
(g_0)$, where $\Phi$ is some linear fractional transformation.
\end{teo}
\begin{proof}
Reasoning as in Theorem \ref{inverse}, there exists a real number $t$ and an
immersion $\phi:\S^2\longrightarrow\H^3$ such that
\begin{itemize}
\item the eigenvalues $\overline{\lambda}_i$ of the new metric $e^{2t}g$ satisfy $\overline{\lambda}_i=e^{-2t}\lambda_i<1/2$ on $\S^2$,
\item the horospherical metric of $\phi$ is $e^{2t}g$,
\item the principal curvatures $k_i$ of $\phi$ satisfy
$e^{-2t}\lambda_i=\frac{1}{2}-\frac{1}{1-k_i}$.
\end{itemize}

Therefore, $\phi$ is a Weingarten immersion with
$$
\widetilde{\cW}(k_1,k_2)=\cW\left(e^{2t}\left(\frac{1}{2}-\frac{1}{1-k_1}\right),
e^{2t}\left(\frac{1}{2}-\frac{1}{1-k_2}\right)\right)=0.
$$

Now, using that
$$
\frac{\parc\widetilde{\cW}}{\parc k_1}\frac{\parc \widetilde{\cW}}{\parc k_2}
=\frac{e^{4t}}{(1-k_1)^2(1-k_2)^2}\frac{\parc\cW}{\parc \landa_1}\frac{\parc
\cW}{\parc \landa_2}
>0
$$
whenever $k_1=k_2$, we have that $\phi$ must be a totally umbilical immersion
\cite{HW,Chern}. Let us remark here that \cite{HW,Chern} only consider surfaces in $\R^3$, but their arguments carry over essentially unchanged to the space forms $\S^3$, $\H^3$, since they rely only on the Codazzi equation, which is the same in all those spaces.

Thus, the result follows easily.
\end{proof}

\bigskip

\noindent {\bf A closing remark:} at the core of the results in the present paper
there is an important phenomenon that we would like to emphasize.

It is well known that the metric geometry of $\H^{n+1}$ is tightly linked to the
conformal geometry of $\S^n$, in the following sense: any isometry of $\H^{n+1}$
induces a conformal transformation of the ideal boundary $\parc_{\8} \H^{n+1} \equiv
\S^n$, and conversely, any conformal transformation of $\S^n$ extends to an isometry
of the unit ball $(\mathbb{B}^{n+1},ds_P^2)\equiv \H^{n+1}$ with respect to the
Poincaré metric.

In this paper, we have shown that this relation is the base of a much more general
situation. Indeed, we have proved that the conformal metrics $g= e^{2\rho} g_0$ on
$\S^n$ are in correspondence with the compact immersed hypersurfaces of $\H^{n+1}$
with prescribed regular hyperbolic Gauss map. Moreover, the fundamental conformal
invariants of $g$ in $\S^n$, i.e. the eigenvalues of the Schouten tensor ${\rm
Sch}_g$, are linked to the basic extrinsic quantities of the hypersurface in
$\H^{n+1}$ (the principal curvatures) by the simple relation $$2\landa_i
=\frac{\k_i+1}{\k_i -1}.$$ It is our impression that this is a far reaching relation
that goes beyond the specific results proved here, which seem to be just the most
visible consequences of the above correspondence.

\def\refname{References}

\end{document}